
\def\LaTeX{%
  \let\Begin\begin
  \let\End\end
  \let\salta\relax
  \let\finqui\relax
  \let\futuro\relax}

\def\UK{\def\our{our}\let\sz s}
\def\USA{\def\our{or}\let\sz z}

\UK



\LaTeX

\USA


\salta

\documentclass[twoside,12pt]{article}
\setlength{\textheight}{24cm}
\setlength{\textwidth}{16cm}
\setlength{\oddsidemargin}{2mm}
\setlength{\evensidemargin}{2mm}
\setlength{\topmargin}{-15mm}
\parskip2mm


\usepackage[usenames,dvipsnames]{color}
\usepackage{amsmath}
\usepackage{amsthm}
\usepackage{amssymb, bbm}
\usepackage[mathcal]{euscript}
\usepackage{hyperref}
%
%


\definecolor{viola}{rgb}{0.3,0,0.7}
\definecolor{ciclamino}{rgb}{0.5,0,0.5}
\definecolor{rosso}{rgb}{0.85,0,0}

\def\rev #1{{\color{red}#1}}

\def\juerg #1{#1}
\def\gianni #1{#1}
\def\an #1{#1}
\def\anold #1{#1}
\def\rev #1{#1}




\bibliographystyle{plain}


%

\finqui

\def\Beq{\Begin{equation}}
\def\Eeq{\End{equation}}
\def\Bsist{\Begin{eqnarray}}
\def\Esist{\End{eqnarray}}

\def\Bthm{\Begin{theorem}}
\def\Ethm{\End{theorem}}
\def\Blem{\Begin{lemma}}
\def\Elem{\End{lemma}}

\def\Bdim{\Begin{proof}}
\def\Edim{\End{proof}}
\def\Bcenter{\Begin{center}}
\def\Ecenter{\End{center}}
\let\non\nonumber




\def\step #1 \par{\medskip\noindent{\bf #1.}\quad}


\def\Lip{Lip\-schitz}
\def\Holder{H\"older}
\def\frechet{Fr\'echet}
\def\aand{\quad\hbox{and}\quad}

\def\lhs{left-hand side}
\def\rhs{right-hand side}

\def\nbh{neighborhood}


\def\nbh{neighb\our hood}


\def\multibold #1{\def\arg{#1}%
  \ifx\arg\pto \let\next\relax
  \else
  \def\next{\expandafter
    \def\csname #1#1#1\endcsname{{\boldsymbol #1}}%
    \multibold}%
  \fi \next}

\def\pto{.}

\def\multical #1{\def\arg{#1}%
  \ifx\arg\pto \let\next\relax
  \else
  \def\next{\expandafter
    \def\csname cal#1\endcsname{{\cal #1}}%
    \multical}%
  \fi \next}

\def\multigrass #1{\def\arg{#1}%
  \ifx\arg\pto \let\next\relax
  \else
  \def\next{\expandafter
    \def\csname gr#1\endcsname{{\mathbb #1}}%
    \multigrass}%
  \fi \next}


\def\multimathop #1 {\def\arg{#1}%
  \ifx\arg\pto \let\next\relax
  \else
  \def\next{\expandafter
    \def\csname #1\endcsname{\mathop{\rm #1}\nolimits}%
    \multimathop}%
  \fi \next}

\multibold
qwertyuiopasdfghjklzxcvbnmQWERTYUIOPASDFGHJKLZXCVBNM.

\multical
QWERTYUIOPASDFGHJKLZXCVBNM.

\multigrass
QWERTYUIOPASDFGHJKLZXCVBNM.

\multimathop
diag dist div dom mean meas sign supp .


\def\accorpa #1#2{\eqref{#1}--\eqref{#2}}
\def\Accorpa #1#2 #3 {\gdef #1{\eqref{#2}--\eqref{#3}}%
  \wlog{}\wlog{\string #1 -> #2 - #3}\wlog{}}


\def\separa{\noalign{\allowbreak}}

\def\neto{\mathrel{{\scriptscriptstyle\nearrow}}}
\def\seto{\mathrel{{\scriptscriptstyle\searrow}}}

\def\graffe #1{\mathopen\{#1\mathclose\}}

\def\<#1>{\mathopen\langle #1\mathclose\rangle}
\def\norma #1{\mathopen \| #1\mathclose \|}

\def\[#1]{\mathopen\langle\!\langle #1\mathclose\rangle\!\rangle}

\def\iot {\int_0^t}
\def\ioT {\int_0^T}
\def\intQt{\int_{Q_t}}
\def\intQ{\int_Q}
\def\iO{\int_\Omega}

\def\dt{\partial_t}
\def\dn{\partial_\nnn}

\def\cpto{\,\cdot\,}

\def\checkmmode #1{\relax\ifmmode\hbox{#1}\else{#1}\fi}
\def\aeO{\checkmmode{a.e.\ in~$\Omega$}}
\def\aeQ{\checkmmode{a.e.\ in~$Q$}}

\def\aeS{\checkmmode{a.e.\ on~$\Sigma$}}
\def\aet{\checkmmode{a.e.\ in~$(0,T)$}}

\def\aat{\checkmmode{for a.a.~$t\in(0,T)$}}


\def\erre{{\mathbb{R}}}




\def\genspazio #1#2#3#4#5{#1^{#2}(#5,#4;#3)}
\def\spazio #1#2#3{\genspazio {#1}{#2}{#3}T0}

\def\L {\spazio L}
\def\H {\spazio H}
\def\W {\spazio W}

\def\C #1#2{C^{#1}([0,T];#2)}


\def\Lx #1{L^{#1}(\Omega)}
\def\Hx #1{H^{#1}(\Omega)}

\def\LQ #1{L^{#1}(Q)}

\def\Luno{\Lx 1}
\def\Ldue{\Lx 2}
\def\Linfty{\Lx\infty}

\def\Huno{\Hx 1}
\def\Hdue{\Hx 2}


\let\badphi\phi
\let\badtheta\theta
\let\theta\vartheta
\let\badeps\epsilon
\let\eps\varepsilon
\let\phi\varphi

\let\TeXchi\chi                         
\newbox\chibox
\setbox0 \hbox{\mathsurround0pt $\TeXchi$}
\setbox\chibox \hbox{\raise\dp0 \box 0 }
\def\chi{\copy\chibox}



\let\emb\hookrightarrow
\def\CO{C_\Omega}
\def\cd{c_\delta}

\def\tarQ{\badphi^Q}
\def\tarO{\badphi^\Omega}
\def\tarsQ{\varsigma^Q}
\def\tarsO{\varsigma^\Omega}

\def\phiz{\phi_0}
\def\sigmaz{\sigma_0}

\def\muO{\mu_\Omega}
\def\etaO{\eta_\Omega}
\def\thetaO{\theta_\Omega}

\def\ustar{u^*}
\def\phistar{\phi^*}
\def\mustar{\mu^*}
\def\sigmastar{\sigma^*}

\def\un{u_n}
\def\phin{\phi_n}
\def\mun{\mu_n}
\def\sigman{\sigma_n}

\def\phih{\hat\phi}
\def\muh{\hat\mu}
\def\sigmah{\hat\sigma}

\def\peps{p^\eps}

\def\reps{r^\eps}
\def\zeps{z^\eps}

\def\soluz{(\phi,\mu,\sigma)}
\def\soluzn{(\phin,\mun,\sigman)}
\def\soluzl{(\psi,\eta,\zeta)}
\def\soluzh{(\phih,\muh,\sigmah)}
\def\soluzstar{(\phistar,\mustar,\sigmastar)}
\def\soluza{(p,q,r)}
\def\soluzaz{(z,p,r)}
\def\soluzeps{(\zeps,\peps,\reps)}

\let\lam\lambda
\def\lamuno{\lam_1}
\def\lamdue{\lam_2}

\def\umin{u_{\min}}
\def\umax{u_{\max}}
\def\Uad{\calU_{ad}}
\def\UR{\calU_R}

\def\Vp{V^*}
\def\calVp{\calV^{\,*}}

\def\normaV #1{\norma{#1}_V}

\let\hat\widehat

\def\mobm{{\mathbbm{m}}}
\def\mobn{{\mathbbm{n}}}

\Begin{document}


%
\title{\juerg{Nutrient control for a viscous\\
Cahn--Hilliard--Keller--Segel model}\\
\juerg{with logistic source describing tumor growth}
}
%
\author{}
\date{}
\maketitle
\Bcenter
\vskip-1cm
{\large\sc Gianni Gilardi$^{(1)}$}\\
{\normalsize e-mail: {\tt gianni.gilardi@unipv.it}}\\[.25cm]
{\large\sc Andrea Signori$^{(2)}$}\\
{\normalsize e-mail: {\tt andrea.signori@polimi.it}}\\[0.25cm]
{\large\sc J\"urgen Sprekels$^{(3)}$}\\
{\normalsize e-mail: {\tt juergen.sprekels@wias-berlin.de}}\\[.5cm]
$^{(1)}$
{\small Dipartimento di Matematica ``F. Casorati'', Universit\`a di Pavia}\\
{\small and Research Associate at the IMATI -- C.N.R. Pavia}\\
{\small via Ferrata 5, I-27100 Pavia, Italy}\\[.3cm] 
$^{(\anold{2})}$
{\small Dipartimento di Matematica, Politecnico di Milano}\\
{\small via Bonardi 9, \anold{I-}20133 Milano, Italy}\\[.2cm] 
$^{(3)}$
{\small Department of Mathematics}\\
{\small Humboldt-Universit\"at zu Berlin}\\
{\small Unter den Linden 6, D-10099 Berlin, Germany}\\
{\small and}\\
{\small Weierstrass Institute for Applied Analysis and Stochastics}\\
{\small Mohrenstrasse 39, D-10117 Berlin, Germany}\\[10mm]
\Ecenter
\begin{center}
\anold{\emph{Dedicated to Pierluigi Colli
on the occasion of his $65^{\rm th}$ birthday,\\ with friendship and admiration}}
\end{center}

\Begin{abstract}
\noindent
In this paper, we address a distributed control problem for a system of partial differential equations describing the evolution of a tumor 
\juerg{that takes the biological mechanism of chemotaxis into account}.
The system describing the evolution is obtained as a nontrivial combination of a Cahn--Hilliard type system accounting for the segregation between tumor cells and healthy cells, with a Keller--Segel type equation accounting for the evolution of a nutrient species and modeling the chemotaxis phenomenon.
First, we \juerg{develop} a robust mathematical background that allows us to analyze an associated optimal control problem. This \juerg{analysis}
forced us to select a source term of logistic type in the nutrient equation and to \juerg{restrict the analysis to the case of} 
two space dimensions.
Then, the existence of \juerg{an} optimal control and \juerg{first-order} necessary conditions for optimality are established.
\vskip3mm
\noindent {\bf Keywords:}
Cahn--Hilliard equation, optimal control, Keller--Segel equation, che\-mo\-taxis\an{, tumor growth.}
\vskip3mm
\noindent {\bf AMS (MOS) Subject Classification:} 
35K55, 
35K61, 
49J20, 
49J50, 
49K20. 
\End{abstract}
\salta
\pagestyle{myheadings}
\newcommand\testopari{\sc Gilardi \ --- \ Signori \ --- \ Sprekels}
\newcommand\testodispari{{\sc \anold{control problem for a Cahn--Hilliard--Keller--Segel model}}}
\markboth{\testodispari}{\testopari}
\finqui
%

\section{Introduction}
\label{Intro}
\setcounter{equation}{0}

Let $\Omega\subset\erre^2$ be some  bounded domain
possessing a smooth boundary $\Gamma:=\partial\Omega$ and the outward unit normal field~$\,\nnn$. 
Denoting by $\dn$ the directional derivative in the direction of~$\nnn$, and putting, with a fixed final time $T>0$,
\Beq
  Q := \Omega\times(0,T) 
  \aand
  \Sigma := \Gamma\times(0,T),
  \non
\Eeq
we study in this paper a {\em distributed control problem\/} associated with \an{(a slight simplification of)} the {\em state system\/} 
given by the following initial-boundary value problem:
\begin{alignat*}{2}
  & \dt\phi - \div\big(\mobm(\phi,\sigma) \nabla \mu\big) 
    =  \gamma(\phi,\sigma) - m \phi && \quad \hbox{in $Q$}\,,\\
  &  \tau \dt \phi - \badeps \Delta \phi + \badeps^{-1}  F'(\phi) - \chi\sigma = \mu && \quad \hbox{in $Q$}\,,\\
  & \dt \sigma -\div \big(\sigma \mobn(\phi,\sigma) \nabla (\ln \sigma + \chi (1-\phi) )\big) 
   =  \beta(\phi) (\kappa_0\sigma - \kappa_\infty \sigma^2) + u&& \quad \hbox{in $Q$}\,,\\
  & \dn \phi = (\mobm(\phi,\sigma)\nabla \mu) \cdot \nnn = (\sigma \mobn(\phi,\sigma) \nabla (\ln \sigma + \chi (1-\phi)) \cdot \nnn = 0
  && \quad \hbox{on $\Sigma$}\,,
  \\
  & \phi(0) = \phiz
  \aand
  \sigma(0) = \sigmaz
  && \quad \hbox{in $\Omega$}\,.
\end{alignat*}
Let us explain the physical meaning of the above symbols.
To begin with, $\phi$ is an order parameter, also referred to as the phase field, which represents
the difference between the \juerg{volume fractions of the tumor and healthy cells}. 
It is normalized in such a way that, at least ideally,
the level sets $\{\phi=1\}:= \{x \in \Omega : \phi (x) = 1\}$ 
and $\{\phi = -1\}$ depict the regions occupied by the pure phases: the tumor and \juerg{healthy tissues},
respectively. 
\juerg{These} regions are then separated by a narrow transition layer
whose thickness scales as the relaxation parameter $\badeps\in(0,1)$.
In the first and third equations, the functions $\mobm(\cdot,\cdot)$ and $\mobn(\cdot,\cdot)$ are nonnegative mobility functions.
As \juerg{it is common for} Cahn--Hilliard type systems, the variable $\mu$ indicates the chemical potential 
associated to the order parameter, $\tau\dt \phi$ is a viscosity contribution, 
and $F'$ denotes the derivative of a configuration potential \juerg{having} a double-well shape.
The mass of the tumor is not conserved, \juerg{which} is captured by the occurrence 
of a source term $S(\phi,\sigma):= \gamma(\phi,\sigma) - m \phi$ on the \rhs\ of the first equation,
where $\gamma$ is a smooth real function on~$\erre^2$, and
$m$ is a positive constant.
Next, the third equation describes the evolution of a nutrient species $\sigma$, where $\chi \geq 0$ denotes the chemotaxis sensitivity, and 
$u:Q \to \erre$ is a control variable. 
The chemotaxis is modeled as in the celebrated Keller--Segel type coupling 
(see, e.g., \cite{H} and the references therein), that is, though the nonlinear term $\chi \div (\sigma \nabla \phi)$ 
\juerg{occurring} in the third equation.
There, $\kappa_0,\kappa_\infty >0$, and $\beta$ is a positive function 
balancing the evolution/saturation effect of the logistic growth and the phase field. 
The choice of a logistic source for the nutrient variable is very common in the Keller--Segel literature: 
see, e.g., \cite{Winkler1, Winkler2, Winkler3} and the references therein, 
as it is a key ingredient to prevent \juerg{a blow-up} of the solution in final time.
Finally, $\phiz$ and $\sigmaz$ are prescribed initial data.


It can be shown (cf.~\cite{RSS}) that the above system is connected to the \juerg{free energy functional}
\begin{align}
\label{freeenergy}
  {\cal F} (\phi, \sigma)
    = \underbrace{\frac \badeps2 \int_\Omega |\nabla \phi |^2 + \frac 1 \badeps\iO F(\phi)}_{=:{\cal E}(\phi)}  
      + \underbrace{\iO \big( \sigma (\ln \an{(\sigma)} - 1) 
      + \chi  \sigma {(1-\phi)} \big) }_{=: {\cal M}(\phi,\sigma)},
\end{align}
where $\cal E$ is the standard Ginzburg--Landau energy approximating the perimeter functional, 
while $\cal M$ is related to the chemotaxis mechanism.
It is worth noticing that, at least formally, the first term of the latter entails a positivity property for~$\sigma$ 
and that the last one, in principle, does not possess a fixed sign, unless one can guarantee that $\phi \in [-1,1]$. 
The last condition, despite \juerg{of} being expected 
from the modeling, may not be fulfilled if the confining potential is defined on the whole real line, 
whereas \juerg{it} directly follows if the potential is singular (cf.~\eqref{regphi}).

A slightly more general version of the above system, uncontrolled and without viscosity, 
i.e., with $\tau =0 $ and $u\equiv 0$, has been addressed in \cite{RSS} from the \juerg{viewpoint of analysis}.
We also refer to \cite{ALL} and \cite{AS}, where a multiphase generalization of the \juerg{ above model}
is introduced using variational principles complying with the second law of thermodynamics in isothermal situations, 
and then analyzed, respectively.
The advantage of the multiphase extension is the possibility of including in the \juerg{model}
further biological effects like angiogenesis 
(see, e.g., \cite{AC, FK}) and \juerg{necrosis}. 
\juerg{In~connection with} the mathematical study of the above system, 
we mention \cite{EG, GARL_2, GLSS,  KS2, SS}, 
where the chemotaxis effects are accounted for through a cross-diffusion type coupling instead.

As they will not play any role from the viewpoint of mathematical investigation, we set for simplicity
$\badeps = \chi = \kappa_0 = \kappa_\infty=1$\an{, and consider the simplified case $\mobm\equiv\mobn \equiv\beta \equiv 1$. 
Thus, the} above system reduces~to
\begin{alignat}{2}
  & \dt\phi - \Delta\mu + m\phi
  = \gamma(\phi,\sigma)
  && \quad \hbox{in $Q$}\,,
  \label{Iprima}
  \\
  & \tau \dt\phi - \Delta\phi + F'(\phi) - \sigma
  = \mu
  && \quad \hbox{in $Q$}\,,
  \label{Iseconda}
  \\
  & \dt\sigma - \Delta\sigma + \div(\sigma\nabla\phi)
  = \sigma - \sigma^2 + u
  && \quad \hbox{in $Q$}\,,
  \label{Iterza}
  \\
  & \dn\phi = \dn\mu = \dn\sigma = 0
  && \quad \hbox{on $\Sigma$}\,,
  \label{Ibc}
  \\
  & \phi(0) = \phiz
  \aand
  \sigma(0) = \sigmaz
  && \quad \hbox{in $\Omega$}\,.
  \label{Icauchy}
\end{alignat}
\Accorpa\Ipbl Iprima Icauchy
\juerg{Although} the well-posedness of the system can be shown for a broad class of potentials, 
including polynomial-type potentials, we will be forced to work under the framework of regular and smooth potentials.
For this reason, we require $F$ to be the Flory--Huggins double-well potential, 
also known as the {\em logarithmic potential}, which is defined as
\begin{align*}
	 F(s) := ((1+s)\ln (1+s)+(1-s)\ln (1-s)) - c_0 s^2 \,,
  \quad s \in (-1,1),
\end{align*}
where $c_0>1$ so that $F$ is nonconvex, with the convention that $0 \ln(0) := \lim_{r \searrow 0} r \ln(r) = 0$. 
Besides, $F$~is extended outside the physical
interval $(-1,1)$ in the usual manner, that is, by continuity at the endpoints $-1$ and~$1$, 
and \juerg{by} $+\infty$ otherwise, to preserve semicontinuity.

To specify  the optimal control problem \juerg{under} study, let us introduce the following tracking-type cost functional:
\begin{align}
  & \calJ(u{;}\phi,\sigma)
  := \frac {\alpha_1} 2 \intQ |\phi-\tarQ|^2
  + \frac {\alpha_2} 2 \iO |\phi(T)-\tarO|^2
  \non
  \\
  & \quad
  + \frac {\alpha_3} 2 \intQ |\sigma-\tarsQ|^2
  + \frac {\alpha_4} 2 \iO |\sigma(T)-\tarsO|^2
  + \frac {\alpha_5} 2 \intQ |u|^2\,,
  \label{Icost}
\end{align}
where the control variable $u$ varies in a proper set of admissible controls (cf.~\eqref{defUad}) 
and $\phi$ and $\sigma$ are the components of the solution $\soluz$ to problem \Ipbl\ 
corresponding to the control~$u$. 
Above, $\alpha_i$, $1\leq i\leq 5$, are nonnegative constants
(not all zero to avoid a trivial situation),
and $\tarQ$, $\tarO$, $\tarsQ$ and $\tarsO$ are given target functions.

Concerning optimal control problems for similar state systems modeling tumor growth, 
we \juerg{refer to} \cite{CSS, CSS2, CSS4, RSS1, S}.
As it will be clarified later on, the restriction \juerg{of} our analysis 
to the two-dimensional case is due to the nonlinear Keller--Segel type coupling, 
which prevents us to infer enough regularity in order to cover the optimal control investigation in three dimensions.
\rev{In particular, even with the help of the viscosity contribution $\tau \dt \phi$ in equation \eqref{Iseconda}, 
in dimension three, the validity of the separation property cannot be established.
This introduces a severe obstruction in the mathematical analysis of the associated optimal control problem.
In fact, the singularity of the double-well potential, connected to a low regularity framework for solutions, 
prevents obtaining robust enough continuous dependence results with respect to the control variable $u$ 
that are a crucial element to deduce the optimality conditions for the minimization problem.}

\section{Statement of the problem and results}
\label{STATEMENT}
\setcounter{equation}{0}

In this section, we state precise assumptions\anold{, set} notations\anold{,} and present our results.
First of all, the set $\Omega\subset\erre^2$ 
is~assumed to be bounded, connected and smooth.
We denote its Lebesgue measure by~$|\Omega|$.
As in the Introduction, $\dn$~stands for the outward normal derivative on $\Gamma:=\partial\Omega$.
Next, if $X$ is a Banach space, then $\norma\cpto_X$ denotes its norm,
with the only exception \juerg{of} the norms in the space $H$ defined below\anold{, whose norm will be indicated by $\norma\cpto$ (i.e., without any subscript),}
and \juerg{in} the $L^p$ spaces ($1\leq p\leq\infty$) constructed on $\Omega$ and~$Q$, which \anold{will be denoted by} $\norma\cpto_p\,$. 
Moreover, in order to simplify notation, the same symbol used for some norm in $X$ 
will also stand for the norm in~$\juerg{X^2:=X\times X}$.
Similarly, if no confusion can arise, we simply write, e.g., $\L2X$ in place of $\L2{X^2}$.
Furthermore, for every Banach space~$X$, the symbols $X^*$ and $\<\cpto,\cpto>_X$ 
denote the dual space of $X$ and the \anold{duality} pairing between $X^*$ and~$X$, respectively.
We also introduce the shorthand
\Beq
  H := \Ldue \,, \quad  
  V := \Huno\,, 
  \aand
  W := \graffe{ v\in\Hdue: \ \dn v=0 }\anold{.}
  \label{defspazi}
\Eeq
Some of our statements \an{involve} the dual space~$\Vp$.
It is understood that we adopt the framework of the Hilbert triplet $(V,H,\Vp)$
obtained by identify\juerg{ing $H$ with} a subspace of $\Vp$ in the usual way\anold{, namely,} in order 
that $\<z,v>_V=\iO \anold{z}v$ for every $z\in H$ and $v\in V$. 

\vskip 2mm

Now, \anold{let us} list
the structur\anold{al assumptions we \juerg{postulate}:}
\begin{align}
  & \hbox{$\tau$ and $m$ are positive real numbers\,.}
  \label{hptm}
  \\
  & \gamma \in C^2(\erre^2;\erre)
  \quad \hbox{is bounded along with its first and second derivatives}
  \non
  \\
  & \quad \hbox{and satisfies} \quad
  \sup|\gamma| < m\,.
  \label{hpgamma}
  \\
  & \hbox{$F:\erre\to(-\infty,+\infty]$ is the logarithmic potential, \anold{that is},}
  \non
  \\
  & \quad F(s) :=
  \left\{
    \begin{array}{ll}
      (1+s)\ln (1+s)+(1-s)\ln (1-s) - c_0 s^2 
      & \quad \hbox{if $|\anold{s}|\leq1$}
      \\
      +\infty
      & \quad \hbox{if $|\anold{s}|>1$}
	\end{array}
  \right.\,,
  \label{defF}
\end{align}
\anold{for} some given real constant $c_0>1$ and the convention \anold{that} $0\ln(0)=0$.
\Accorpa\HPstruttura hptm defF

For the data, we make the assumptions listed below.
Even though the control $u$ has to be consider\anold{ed as} a fixed datum at the present stage,
it is convenient to introduce a constant $M$ which is an upper bound for its $L^\infty$ norm.
Thus, we \anold{require}~that:
\begin{align}
  & \hbox{$M$ is a positive constant, and $u\in\LQ\infty$ satisfies}
  \,\,
  \norma u_\infty \leq M\,.
  \label{hpu}
  \\
  & \phiz \in \Hx3 \cap W
  \aand
  \norma\phiz_\infty < 1 \,.
  \label{hpphiz}
  \\
  & \sigmaz \in V \cap \Linfty
  \aand
  \sigmaz \geq 0 \quad \aeO\,.
  \label{hpsigmaz}
\end{align}
\Accorpa\HPdati hpu hpsigmaz

The above assumptions guarantee that the problem \Ipbl\ is well\anold{-}posed 
in a proper functional framework and that satisfactory stability and continuous dependence results hold true.
However, for convenience, we consider the variational formulation of the problem \anold{that is presented below}.

\Bthm
\label{Wellposedness}
Assume \HPstruttura\ on the structure of the system and \HPdati\ on the data.
Then there exists a unique triplet $\soluz$ enjoying the properties
\begin{align}
  & \phi \in \W{1,\infty}{\rev V} \cap \H1W \cap \L\infty{\Hx{\rev3}}
  \aand
  \norma\phi_\infty < 1\,,
  \label{regphi}
  \\
  & \mu \in \L\infty W\,, 
  \label{regmu}
  \\
  & \sigma \in \H1H \cap \C0V \cap \L2W \cap \LQ\infty
  \aand
  \sigma \geq 0 \quad \aeQ\,,
  \label{regsigma}
\end{align}
\Accorpa\Regsoluz regphi separation
and solving the variational problem
\begin{align}
  & \iO \dt \phi \, v
  + \iO \nabla\mu \cdot \nabla v
  + m \iO \phi v
  = \iO \gamma(\phi,\sigma) \, v
  \non
  \\
  & \quad \hbox{for every $v\in V$ and \aet}\,,
  \label{prima}
  \\
  & \tau \iO \dt\phi \, v
  + \anold{\iO} \nabla\phi \cdot \nabla v
  + \iO F'(\phi) v
  - \iO \sigma v
  = \iO \mu v
  \non
  \\
  & 
  \quad \hbox{for every $v\in V$ and \aet}\,,
  \label{seconda}
  \\
  & \iO \dt\sigma \, v
  + \iO \nabla\sigma \cdot \nabla v
  - \iO \sigma \nabla\phi \cdot \nabla v
  = \iO \bigl( \sigma - \sigma^2 + u \bigr) v
  \non
  \\
  & 
  \quad \hbox{for every $v\in V$ and \aet}\,,
  \label{terza}
  \\
  & \phi(0) = \phiz
  \aand
  \sigma(0) = \sigmaz \,.
  \label{cauchy}
\end{align}
\Accorpa\Pbl prima cauchy
In particular, the solution satisfies the initial-boundary value problem \Ipbl,
the equations and the boundary conditions holding \aeQ\ and \aeS, respectively.
Moreover, \anold{it} satisfies the stability estimate
\begin{align}
  & \norma\phi_{\W{1,\infty}H\cap\H1W\cap\L\infty{\anold{\Hx{\rev3}}}}
  + \norma\mu_{\L\infty W}
  \non
  \\
  & \quad
  + \norma\sigma_{\H1H\cap\C0V\cap\L2W\cap\LQ\infty}
  \leq K_1\,,
  \label{stab}
\end{align}
\anold{as well as the separation property
\begin{align}
  |\phi| \leq 1-\delta_0
  \quad \aeQ\,,
  \label{separation}
\end{align}
with} some constants $K_1>0$ and $\delta_0\in(0,1)$ that depend only on the structure of the system,
$\Omega$, $T$, the initial data, and the constant~$M$.
In particular, they are independent of~$u$.
\Ethm

\Bthm
\label{Contdep}
Under the assumptions of Theorem~\ref{Wellposedness} on the structure and the initial data,
let $u_i\in\LQ\infty$, $i=1,2$, \anold{be such that} $\norma{u_i}_{\anold{\infty}}\leq M$, and let
$(\phi_i,\mu_i,\sigma_i)$ be the corresponding solutions.
Then the inequality
\begin{align}
  & \norma{\phi_1-\phi_2}_{\H1H\cap\L\infty V\cap\L2W}
  + \norma{\mu_1-\mu_2}_{\L2W}
  \non
  \\
  & \quad
  + \norma{\sigma_1-\sigma_2}_{\L\infty H\cap\L2V}
  \leq K_2 \norma{u_1-u_2}_{\L2H}
  \label{contdep}
\end{align}
holds true with a positive constant $K_2$ that depends only on the structure of the system,
$\Omega$, $T$, the initial data, and the constant~$M$.
\Ethm

By accounting for the above results, we deal with the control problem sketched in the Introduction
and related to the cost functional
\begin{align}
  & \calJ(u\anold{;}\phi,\sigma)
  := \frac {\alpha_1} 2 \intQ |\phi-\tarQ|^2
  + \frac {\alpha_2} 2 \iO |\phi(T)-\tarO|^2
  \non
  \\
  & \quad
  + \frac {\alpha_3} 2 \intQ |\sigma-\tarsQ|^2
  + \frac {\alpha_4} 2 \iO |\sigma(T)-\tarsO|^2
  + \frac {\alpha_5} 2 \intQ |u|^2\,,
  \label{cost}
\end{align}
where we assume that
\begin{align}
  & \hbox{$\alpha_1,\dots,\alpha_5$ are nonnegative constants\,,}
  \label{hpalpha}
  \\
  & \tarQ, \, \tarsQ \in \LQ2
  \aand
  \tarO, \, \tarsO \in \Ldue \,.
  \label{hptarget}
\end{align}
For the distributed control variable~$u$, we choose as set of {\anold{\it admissible controls}}
\Beq
  \Uad : = \big\{ u \in {\cal U}: \umin \leq u \leq \umax \quad \aeQ \big\}\,,
  \quad \hbox{where} \quad 
  \calU := \LQ\infty\,,
  \label{defUad}
\Eeq
\anold{and} assum\anold{e} that 
\Beq
  \umin, \, \umax \in \LQ\infty
  \quad\mbox{with}\quad
  \umin \leq \umax \quad \aeQ \,.
  \label{hpuminmax}
\Eeq
Notice that $\Uad$ is nonempty and bounded in $\LQ\infty$.
More precisely, our well-posedness and continuous dependence results
can be applied by taking $M$ as, e.g., the maximum between $\norma\umin_\infty$ and $\norma\umax_\infty$.
\anold{Summing up}, the control problem \juerg{under study} reads as follows:
\begin{align}
  & \hbox{\sl Minimize the cost functional \eqref{cost} subject to $u\in\Uad$ and}
  \non
  \\
  & {}\quad \hbox{\sl the solution  $\soluz$ to \Pbl\ corresponding to $u$\anold{.}}
  \label{ctrl}
\end{align}
This problem is \anold{analyzed} in detail in Section~\ref{CONTROL}, 
where we first prove the existence of an optimal strategy
and then provide a \anold{first-order} necessary condition for an element $\ustar\in\Uad$ to be an optimal control.
The \anold{latter is encoded in} the variational inequality
\Beq
  \intQ (r+\alpha_5\ustar) (u-\ustar) \geq 0
  \quad \hbox{for every $u\in\Uad$}\,,
  \non
\Eeq
where $r$ is the third component of the solution $\soluza$ to a \anold{suitable adjoint} problem \anold{to \Ipbl}.
\juerg{We do not  describe it here in detail, since it should be formulated in a suitable weak form whose introduction requires some space.
Instead, we provide a formal \anold{(strong)} version thereof: it is given by} the backward-in-time \anold{parabolic} problem 
\begin{alignat}{2}
  & - \dt(p+\tau q) 
  - \Delta q 
  + (m-\lamuno)p + \lam q
  + \div(\sigmastar\nabla r)
  = g_1 
  && \anold{\quad \hbox{in $Q$}}\,,
  \non
  \\
  & - \Delta p 
  = q
  && \anold{\quad \hbox{in $Q$}}\,,
  \non
  \\
  & - \dt r
  - \Delta r 
  - \lamdue p - q + (2\sigmastar-1) r
  - \nabla\phistar \cdot \nabla r 
  = g_2
  && \anold{\quad \hbox{in $Q$}}\,,
  \non
  \\
  & \dn p = \dn q = \dn r = 0
  && \quad \anold{\hbox{on $\Sigma$}}\,,
  \non
  \\
  & (p+\tau q)(T) = g_3
  \aand
  r(T) = g_4
   && \quad \anold{\hbox{in $\Omega$\,.}}
 \non
\end{alignat}
Here, we \juerg{have used the denotations}
\begin{align}
  & \hbox{$\phistar$ and $\sigmastar$ are the components of the state \anold{system} corresponding to~$\ustar$}\,,
  \non
  \\
  & \lamuno := \gamma_\phi(\phistar,\sigmastar) , \quad
  \lamdue := \gamma_\sigma(\phistar,\sigmastar),
  \aand
  \lam := F''(\phistar) ,
  \non
  \\
  & g_1 := \alpha_1(\phistar-\tarQ) , \quad
  g_2 := \alpha_2(\sigmastar-\tarsQ) , \quad
  g_3 := \alpha_3 \bigl( \phistar(T)-\tarO \bigr),
  \non
  \\
  & \aand
  g_4 := \alpha_4 \bigl( \sigmastar(T)-\tarsO \bigr),
  \non
\end{align}
\anold{where $\gamma_\phi$ and $\gamma_\sigma$ indicate the partial derivatives of $\gamma$ with respect to the respective variables.}

\medskip

The rest of the paper is organized as follows.
We continue the present section by recalling some tools
and stating a general rule on the notation regarding the constants 
that appear in the estimates \juerg{below}.
The next section is devoted to the proofs of Theorems~\ref{Wellposedness}
and~\ref{Contdep} \anold{on the state system},
while the control problem is \anold{discussed} in the last section.

Throughout the paper, we will repeatedly use the Young inequality
\Beq
  a\,b \leq \delta\,a^2 + \frac 1{4\delta} \, b^2
  \quad \hbox{for all $a,b\in\erre$ and $\delta>0$}\,,
  \label{young}
\Eeq
as well as the \Holder, Schwarz and Poincar\'e inequalities.
We recall the latter:
there exists a constant $\CO$\anold{,} depending only on~$\Omega$\anold{,} such~that
\Beq
  \norma v^2
  \leq \CO \bigl( \norma{\nabla v}^2 + |v_\Omega|^2 \bigr)
  \quad \hbox{for every $v\in V$}\,,
  \label{poincare}
\Eeq
where $v_\Omega$ denotes the mean value of~$v$,~i.e.,
\Beq
  v_\Omega := \frac 1 {|\Omega|} \iO v
  \quad \hbox{for $v\in\Luno$}.
  \label{defmean}
\Eeq
\anold{T}he same symbol $v_\Omega$ \anold{will be used} also if $v$ is time dependent.
Moreover, we take advantage of the (two-dimensional) continuous embeddings
\Beq
  V \emb \Lx p
  \quad \hbox{for $p\in[1,+\infty)$}
  \aand
  \L\infty H \cap \L2V \emb \LQ4 \,.
  \label{emb}
\Eeq
Finally, we will employ the abbreviations
\Beq
  Q_t := \Omega\times(0,t)
  \quad \hbox{for $t\in(0,T]$}
  \aand
  Q^t := \Omega\times(t,T)
  \quad \hbox{for $t\in[0,T)$}.
  \label{defQt}
\Eeq

We conclude this section by stating a general rule 
concerning the constants that appear in the estimates we perform in the following\anold{: w}e use the small-case symbol $\,c\,$ for a generic constant
whose actual values may change from line to line and even within the same line
and \juerg{depends} only on~$\Omega$, the structure of the system,
and the constants and the norms of the functions involved in the assumptions of the statements.
In particular, the values of $\,c\,$ may depend on the constant $M$ that appears in \eqref{hpu},
but they are independent of \anold{the control variable} $u$.
A~small-case symbol with a subscript like $\cd$
indicates that the constant may depend on the parameter~$\delta$, in addition.
On the contrary, we mark precise constants that we can refer~to
by using different symbols 
\anold{(e.g.,~a capital letter as in \eqref{poincare}).}


\section{The state system}
\label{WP}
\setcounter{equation}{0}

\subsection{Existence}
\label{EXISTENCE}

This subsection regards the existence part of Theorem~\ref{Wellposedness}
and the stability estimate~\eqref{stab}.
The \anold{nonviscous and uncontrolled} version of our problem\anold{,} i.e.,
\anold{system \Ipbl\ with} $\tau=0$ in \eqref{seconda} \anold{and} $u=0$\anold{,}
is~a particular case of the system studied in~\cite{RSS},
where well-posedness and regularity results have been \anold{established}.
As for existence, the first observation made in the quoted paper regards equation~\eqref{terza},
whose principal part can be written in a different form
by accounting for the identity \anold{(see also \an{the associated free energy} \eqref{freeenergy})}
\Beq
 \nabla\sigma - \sigma\nabla\phi
 = \sigma \nabla(\ln(\sigma) + 1 - \phi)\,,
 \label{identity}
\Eeq
provided that $\sigma$ is assumed to be positive.
Then, formally testing the new version of \eqref{terza} by $\ln(\sigma)+1-\phi$,
one \juerg{obtains} an energy estimate involving~$\ln(\sigma)$,
which implies the positivity of~$\sigma$.
However, this is \juerg{only} formal, and the \anold{a}uthors of \cite{RSS}
proceeded rigorously by performing the energy estimate
on the solution to a suitably regularized problem.
Then, further estimates led to existence and regularity
in a functional \juerg{analytic} framework that is close to the one related to 
our regularity requirements \anold{\eqref{regphi}--\eqref{regsigma} and \eqref{separation}} 
(which could \anold{actually} be improved as in~\cite{RSS}).
\anold{On the other hand}, in the first estimate \anold{of \cite{RSS},} the regularity of the time derivative of~$\phi$ 
\anold{has to be understood as a} $\Vp$-valued function, since $\tau=0$ in that case.
Nevertheless, the whole procedure developed there can be repeated here:
\anold{it is worth noting that}
the additional term $u$ in \eqref{terza} does not cause any trouble, \anold{as we are assuming it to be bounded.
Moreover, the} compatibility conditions on $m$ and $\gamma$ given in \eqref{hpgamma}
can be used in the same way to control the mean value of~$\phi$,
and a better estimate of $\dt\phi$ is obtained from the very beginning
thanks to the presence of the viscosity term $\tau\dt\phi$ in~\eqref{seconda}.
In particular, whenever it is convenient
(e.g., when considering the nonlinear elliptic operator $-\Delta+F'$),
one can move the term $\dt\phi$ to the \rhs\ in our case.
\anold{That said,} some \anold{minor} modifications are needed,
and we have to show how the argument of each step of \cite{RSS} can be adapted to the new situation.
For the sake of brevity, we confine ourselves to present the list of the \anold{modifications} acting on problem \Pbl, directly,
and we just focus on the contributions due to the new term
that \juerg{originate from} formally testing by some functions and integrating with respect to time.

By taking $v=\dt\phi$ in \eqref{seconda}
(coupled with a suitable choice of the test functions in the other equations),
one \anold{obtains} a nonnegative term on the \lhs\
that leads~to
\Beq
  \norma{\dt\phi}_{\L2H} \leq c \,,
  \non
\Eeq
in addition.
\rev{%
The next step relies on differentiating equations \eqref{prima} and \eqref{seconda} with respect to time and test the resulting equalities  by~$\dt\phi$
and $- \Delta \dt \phi$, respectively.
Adding those leads to 
\begin{align*}
	& 
	\frac 12 \frac {d}{dt} \iO|{\dt \phi}|^2
	+\frac {\tau}2 \frac {d}{dt} \iO|{\nabla \dt \phi}|^2
	+ \iO|{\Delta \dt \phi}|^2
	+ m \iO|{\dt \phi}|^2
	= \iO F''(\phi) \dt \phi \Delta \dt \phi 
	\\ & \quad 
	+ \iO \dt \sigma \Delta \dt \phi 
	+ \iO \big(\gamma_\phi(\phi,\sigma) \dt \phi + \gamma_\sigma(\phi,\sigma) \dt \sigma\big) \Delta \dt \phi \,.
\end{align*}
Here, $\gamma_\phi=\gamma_\phi(\phi,\sigma) $ and $\gamma_\sigma=\gamma_\sigma(\phi,\sigma)$ denote the partial derivatives of $\gamma$ with respect to $\phi$ and $\sigma$, respectively.
Notice that, due to \eqref{hpgamma} they both are uniformly bounded.
All the terms on the \rhs\ can be readily controlled as in \cite{RSS} and the difference is that we now have the additional viscosity term 
$\frac {\tau}2 \frac {d}{dt} \iO|{\nabla \dt \phi}|^2$ on the \lhs.
Provided that we can estimate $\dt\phi(0)$ in~$V$, this would produce,
after integration over time, using the Gronwall's lemma and elliptic regularity, that 
\Beq
	\norma{\phi}_{\W{1,\infty} V \cap \H1 W}
	\leq c\,.
	\non
\Eeq
Next, we test \eqref{prima} by $\dt\mu$ and observe that the quantity $\nabla\mu(0)$ has to be controlled as well.
}%
To overcome these difficulties, we subtract \eqref{seconda} \juerg{from} \eqref{prima} multiplied by $\tau$ 
and take $t=0$ in the \juerg{resulting} equality.
Thanks to a cancellation, this leads to the elliptic problem
\begin{align}
  & \tau \iO \nabla\mu(0) \cdot \nabla v
  + \iO \mu(0) \, v 
  \non
  \\
  & = \iO \bigl( \anold{\tau \gamma(\phiz,\sigmaz)- \tau m \phiz - \Delta\phiz   }+ F'(\phiz) - \sigmaz \bigr) \, v
  \quad \hbox{for every $v\in V$}.
  \non
\end{align}
Thus, $\mu(0)$ is bounded in $\Hx3$, by virtue of the regularity theory of elliptic equations
and our assumptions on \juerg{the} structure of the system and the initial data,
and $\dt\phi(0)$ is bounded in $V$ as a consequence
\juerg{of} \eqref{prima} written at $t=0$.
This concludes the formal proof of the existence of a solution, 
\anold{as the rest of the details can be filled in by arguing along the same line of arguments \juerg{as in}~\cite{RSS}}.
Moreover, a~clever inspection of the above argument 
also shows that both the stability estimate \eqref{stab} and the separation property \eqref{separation}
hold true with constants $K_1$ and $\delta_0$ that have the dependences specified in the statement.

\subsection{Uniqueness and continuous dependence}
\label{CONTDEP}

\anold{Next, we move to proving} the uniqueness part of Theorem~\ref{Wellposedness} 
and the continuous dependence \an{presented in} Theorem~\ref{Contdep}.
More precisely, given $u_i$, $i=1,2$, as in \anold{Theorem~\ref{Contdep}},
we first prove \gianni{a continuous dependence estimate (in~the direction of~\eqref{contdep})
by assuming that $(\phi_i,\mu_i,\sigma_i)$ are arbitrary solutions corresponding to~$u_i$.
Then, we derive the uniqueness of the solution as a consequence.
Finally, we complete the proof of~\eqref{contdep}}.
\anold{Along with the fixed} controls and the \anold{corresponding} solutions\anold{, we} set for convenience
\Beq
  u := u_1-u_2 \,, \quad
  \phi := \phi_1-\phi_2 \,, \quad
  \mu := \mu_1-\mu_2\,, 
  \aand
  \sigma := \sigma_1-\sigma_2\,, 
  \non
\Eeq
and notice that these functions satisfy
\begin{alignat}{2}
  & \dt \phi 
  - \Delta\mu 
  + m \phi 
  = \gamma(\phi_1,\sigma_1) - \gamma(\phi_2,\sigma_2) 
  && \quad \anold{\hbox{in $Q$}}\,,
  \label{diffprima}
  \\
  & \tau \dt\phi 
  - \Delta\phi 
  + F'(\phi_1) - F'(\phi_2)
  - \sigma
  = \mu
  && \quad \anold{\hbox{in $Q$}}\,,
  \label{diffseconda}
  \\
  & \dt\sigma 
  - \Delta\sigma 
  + \div(\sigma_1 \nabla\phi_1 - \sigma_2 \nabla\phi_2)
  = \sigma - (\sigma_1^2 - \sigma_2^2) + u 
  && \quad \anold{\hbox{in $Q$}\,.}
  \label{diffterza}
\end{alignat}
Moreover, $\phi$, $\mu$\anold{,} and $\sigma$ satisfy homogeneous Neumann boundary conditions,
and $\phi$ and $\sigma$ vanish at $t=0$.
We multiply the above equations by $\mu$, $\dt\phi{-}\Delta\phi$, and~$\sigma$, respectively,
integrate over~$\Omega$ and by parts in space, sum up, and rearrange.
We also add to both sides the same quantity
$\frac12\frac d{dt}\iO|\phi|^2=\iO\phi\dt\phi$.
Due to a cancellation, we obtain~that
\begin{align}
  & \iO |\nabla\mu|^2
  + \tau \iO |\dt\phi|^2
  + \frac 12 \, \frac d{dt} \, \normaV\phi^2
  + \frac \tau 2 \, \frac d{dt} \iO |\nabla\phi|^2
  + \iO |\Delta\phi|^2
  \non
  \\
  & \quad {}
  + \frac 12 \, \frac d{dt} \iO |\sigma|^2
  + \iO |\nabla\sigma|^2
  \non
  \\
  & = \iO \bigl( -m \phi + \gamma(\phi_1,\sigma_1) - \gamma(\phi_2,\sigma_2) \bigr) \mu
  - \iO \bigl( F'(\phi_1) - F'(\phi_2) \bigr) \bigl( \dt\phi - \Delta\phi \bigr)
  \non
  \\
  & \quad {}
  + \iO \sigma (\dt\phi - \Delta\phi)
  - \iO \mu \Delta\phi
  \anold{{}+ \iO \bigl( \sigma \nabla\phi_1 + \sigma_2 \nabla\phi \bigr) \cdot \nabla\sigma}
  \non
  \\
  & \quad {}
  + \iO \sigma^2
  - \iO (\sigma_1+\sigma_2) |\sigma|^2
  + \iO u \sigma 
  + \iO \phi \, \dt\phi \,.
  \label{testdiff}
\end{align}
In estimating the \rhs\ of \eqref{testdiff},
we \anold{employ} the convention on the constants announced at the end of Section~\ref{STATEMENT}.
\anold{From now on, we} allow the values of the constant $c$ to \juerg{additionally} depend on the \juerg{fixed}
solutions $(\phi_i,\mu_i,\sigma_i)$\anold{, $i=1,2$}.
Later on, we will show how this further dependence can be removed.
We repeatedly make use of the Young and Poincar\'e inequalities,
and account for the \Lip\ continuity of $\gamma$ and~$F'$,
the latter \anold{being fulfilled} on every compact subset of~$(-1,1)$, \anold{due to \eqref{separation}}.
In the sequel, $\delta$~\anold{indicates} an arbitrary positive number \anold{whose value is yet to be selected}.
We have~that
\begin{align*}
  & \iO \bigl( -m \phi + \gamma(\phi_1,\sigma_1) - \gamma(\phi_2,\sigma_2) \bigr) \mu
  \\
  & \anold{=\iO \bigl( -m \phi + \gamma(\phi_1,\sigma_1) - \gamma(\phi_2,\sigma_2) \bigr) (\mu- \mu_\Omega) + \iO \bigl( -m \phi + \gamma(\phi_1,\sigma_1) - \gamma(\phi_2,\sigma_2) \bigr) \mu_\Omega }
  \\
  & \leq \delta \,|\muO|^2
  + \delta \iO |\nabla\mu|^2
  + \cd \iO \bigl( |\phi|^2 + |\sigma|^2 \bigr)\,, 
  \non
\end{align*}
as well as
\Beq
  - \iO \bigl( F'(\phi_1) - F'(\phi_2) \bigr) \bigl( \dt\phi - \Delta\phi \bigr)
  \leq \delta \iO |\dt\phi|^2
  + \delta \iO |\Delta\phi|^2
  + \cd \iO |\phi|^2 .
  \non
\Eeq
Next, \anold{using integration by parts, we find that}
\Beq
  - \iO \mu \Delta\phi
  = \iO \nabla\mu \cdot \nabla\phi
  \leq \delta \iO|\nabla\mu|^2 
  + \cd \iO |\nabla\phi|^2\,,
  \non
\Eeq
and the remaining terms can be easily \anold{controlled. In particular,}
\anold{let us recall} that $\nabla\phi_1$ and $\sigma_2$ are bounded,
and that both $\sigma_1$ and $\sigma_2$ are nonnegative, \anold{as a consequence of Theorem~\ref{Wellposedness}}.
\anold{Then,} we \anold{are left with handling} the mean value $\muO$ of~$\mu$ that occurs on the \anold{\rhs\ of the first estimate above}.
To this end,
we integrate \eqref{diffseconda} over~$\Omega$ and obtain~that
\Beq
  |\Omega| \, \muO
  = \tau \iO \dt\phi
  + \iO \bigl( F'(\phi_1) - F'(\phi_2) \bigr)
  - \iO \sigma \,.
  \non
\Eeq
By squaring and recalling that $\norma v_1^2\leq|\Omega|\,\norma v^2$ for every $v\in H$,
we infer~that
\begin{align}
  & |\Omega|^2 \, |\muO|^2
  \leq 3 \Bigr\{
    \Bigl( \tau \iO |\dt\phi| \Bigr)^2
    + \Bigl( \iO |F'(\phi_1) - F'(\phi_2)| \Bigr)^2
    + \Bigl( \iO \anold{\sigma} \Bigr)^2
  \Bigr\}
  \non
  \\
  & \leq 3 |\Omega| \Bigl(
    \tau^2 \iO |\dt\phi|^2
    + \iO |F'(\phi_1) - F'(\phi_2)|^2
    + \iO |\sigma|^2
  \Bigr)\an{\,,}
  \non
\end{align}
whence\anold{, owing to the \Lip\ continuity of $F'$,} also
\Beq
  \frac {|\Omega|} {6\tau} \, |\muO|^2
  \leq \frac \tau 2 \iO |\dt\phi|^2
  + c \iO |\phi|^2
  + c \iO |\sigma|^2 \,.
  \label{mean}
\Eeq
At this point, we add \eqref{mean} to \eqref{testdiff},
account for the above estimates, and choose $\delta>0$ small enough.
By integrating with respect to time, applying the Gronwall lemma,
and using the Poincar\'e inequality,
we obtain that
\Beq
  \norma\phi_{\H1H\cap\L\infty V}
  + \norma\mu_{\L2V}
  + \norma\sigma_{\L\infty H\cap\L2V}
  \leq c \, \norma u_{\L2H}\anold{,}
  \label{less}
\Eeq
\anold{where} the \anold{appearing} constant $c$ may also depend on the \anold{fixed solutions}.
Nevertheless, from \eqref{less} we can derive uniqueness since the solutions \anold{were} arbitrary.
Therefore, coming back to the proof just performed
and realizing that the solutions must coincide with the ones constructed in the previous subsection,
we can replace the norms of the solutions that enter the estimates
(e.g., $\norma{\sigma_2}_\infty$)
by~owing to the stability estimate \eqref{stab}
and the constant $c$ appearing in \eqref{less} by a constant independent of the solutions,
thus with the dependence specified in the statement of Theorem~\ref{Contdep} for the constant~$K_2$.
Nevertheless, the new version of \eqref{less} is still weaker than~\eqref{contdep}.
\anold{Besides, by a straightforward comparison argument}, we can easily complete the estimate, since \eqref{diffprima}, \eqref{diffseconda}
and what we have already obtained imply~that
\Beq
  \norma{\Delta\phi}_{\L2H} + \norma{\Delta\mu}_{\L2H}
  \leq c \, \norma u_{\L2H} \,.
  \non
\Eeq
Then, \eqref{contdep} follows \juerg{from the} elliptic regularity \anold{theory}.


\section{The control problem}
\label{CONTROL}
\setcounter{equation}{0}

In this section, we \anold{address} the control problem \eqref{ctrl}.
It is understood that all of the assumptions on the structure of the original system, the data,
and the ingredients of the cost functional \eqref{cost}, which we have made throughout the paper, are in force \anold{from now on}.


\subsection{Existence of an optimal strategy}
\label{EXISTCONTROL}

The first result of ours is the following:

\Bthm
\label{OKcontrol}
The \anold{optimization} problem \eqref{ctrl} \anold{admits} at least one solution $\ustar$.
\Ethm
 
\Bdim
We use the direct method of calculus of variations.
\anold{To begin with, let us notice that $\cal J$ is bounded from below by zero.}
Now, we pick a minimizing sequence $\graffe{\un}$ in $\Uad$
and the corresponding sequence $\graffe{\soluzn}$ of solutions to the state system.
Since $\Uad$ is bounded in~$\calU$, we can assume~that
\Beq
  \un \to \ustar
  \quad \hbox{weakly star in $\calU$}
  \non
\Eeq
as $n\neto\infty$, for some limit function $\ustar$, which must belong \anold{to~$\Uad$ since $\Uad$} is convex and \anold{(strongly)} closed.
Besides, the corresponding solutions are bounded as well 
in the topologies specified in Theorem~\ref{Wellposedness}.
\juerg{Therefore it follows, possibly only on a subsequence which is again indexed by $n$,} ~that\anold{, as $n\neto\infty$,}
\Bsist
  && \phin \to \phistar
  \quad \hbox{weakly star in $\W{1,\infty}{\rev V}\cap\H1W\cap\L\infty{\anold{\Hx{\rev3}}$}}\,,
  \non
  \\
  && \mun \to \mustar
  \quad \hbox{\rev{weakly star in $\L\infty W$}}\,,
  \non
  \\
  && \non
  \sigman \to \sigmastar
  \quad \hbox{\anold{weakly star} in $\H1H\cap\L\infty V\cap\L2W\cap\LQ\infty$}\,,
  \qquad
\Esist
with \anold{suitable} limit functions \anold{$\phi^*, \mu^*$, and $\sigma^*$}.
Moreover, the separation property $\norma\phin_\infty\leq1-\delta_0$
is satisfied with some $\delta_0\in(0,1)$ independent of~$n$.
Therefore, it is immediately seen that the triplet $\soluzstar$ is the solution to the state system
corresponding to $\ustar$ and~that
\Beq
  \calJ(\un;\phin,\sigman) \to \calJ(\ustar;\phistar,\sigmastar) \,.
  \non
\Eeq
On the other hand, \anold{by construction,} we also have that
$\lim_{n\neto\infty} \calJ(\un;\phin,\sigman)$
coincides with the infimum of~$\calJ$
since the sequence $\graffe{\un}$ is minimizing~$\calJ$.
Therefore, \anold{the infimum is attained and} $\ustar$ is an optimal control.
\Edim
 

\subsection{The control-to-state mapping}
\label{CTS}

In this section, we introduce the {\it control-to-state} mapping $\calS$\anold{, also referred to as the {\it solution operator},} and prove its \frechet\ differentiability \anold{in a suitable mathematical framework}.
Along with the space $\calU$ and the set $\Uad$ of the admissible controls defined in~\eqref{defUad}, 
we introduce the state space $\calY$ and the open \nbh\ $\UR$ of $\Uad$ by setting
\juerg{
\begin{align}
  & \calY := \calY_1 \times \calY_2 \times \calY_3\,,
  \quad \hbox{where} 
	\quad 
  \calY_1 := \H1H\cap\L\infty V , \quad
  \calY_2 := \L2V\,,
  \non
  \\
  & \aand
  \calY_3 := \anold{\H1 \Vp \cap}\L\infty H\cap\L2V\,,
  \label{defY}
  \\
  & \UR := \graffe{u\in\calU : \ \norma u_\infty < R}\,,
  \quad \hbox{where} \quad
  R := \max \graffe{\norma\umin_\infty\,,\norma\umax_\infty} + 1 \,.
  \label{defUR}
\end{align}
Finally, we define the map
\begin{align}
  & \calS : \UR \to \calY ;  \quad 
  u \mapsto \hbox{the solution $(\phi,\mu,\sigma)$ to \Pbl\ corresponding to $u$}\,.
  \label{defS}
\end{align}
We} notice that we can apply Theorems~\ref{Wellposedness} and \ref{Contdep}
with $M=R$ to ensure that $\calS$ is well defined and obtain the stability estimate \eqref{stab},
the separation property \eqref{separation},
and the continuous dependence estimate \eqref{contdep},
with fixed constants $K_1$, $\delta_0$ and $K_2$ independent of $u\in\UR$.
As a consequence, there exists a positive constant $K_3$ such~that
\begin{align}
  & \norma{F^{(i)}(\phi)}_\infty
  \leq K_3\,,
  \quad \hbox{for $0\leq i\leq 3$ and every $u\in\UR$},
  \non
  \\
  & \quad \hbox{where $\phi$ is the first component of $\calS(u)$}.
  \label{derivF}
\end{align}
The main result of this section is the \frechet\ differentiability of~$\calS$.
This is related to the linearized system \juerg{introduced now. To this end, let} $u\in\UR$ and $\soluz:=\calS(u)$.
We \juerg{denote}, for brevity,
\Beq
  \lamuno := \gamma_\phi(\phi,\sigma) , \quad
  \lamdue := \gamma_\sigma(\phi,\sigma),
  \aand
  \lam := F''(\phi).
  \label{deflambda}
\Eeq 
Then\an{,} the linearized system corresponding to $u$ and to the variation $h\in\calU$
is the system, 
whose unknown is the triplet $\soluzl\in\calY$,
\begin{align}
  & \iO \dt\psi \, v
  + \iO \nabla\eta \cdot \nabla v
  + m \iO \psi v
  \non
  \\
  & = \iO \bigl( \lamuno \psi + \lamdue \zeta \bigr) v
  \quad \hbox{\aet\ and for every $v\in V$}\,,
  \qquad
  \label{primal}
  \\
  & \tau \iO \dt\psi \, v
  + \iO \nabla\psi \cdot \nabla v
  + \iO \bigl( \lam \psi - \zeta \bigr) v
  \non
  \\
  & = \iO \eta v
  \quad \hbox{\aet\ and for every $v\in V$}\,,
  \label{secondal}
  \\
  & 
	\anold{\<\dt\zeta , v>_V}
  + \iO \nabla\zeta \cdot \nabla v
  - \iO \bigl( \zeta\nabla\phi + \sigma\nabla\psi \bigr) \cdot \nabla v
  \non
  \\
  & = \iO \bigl( \zeta - 2\sigma\zeta + h \bigr) v
  \quad \hbox{\aet\ and for every $v\in V$}\,,
  \label{terzal}
  \\
  & \psi(0) = 0
  \aand
  \zeta(0) = 0 \,.
  \label{cauchyl}  
\end{align}
\Accorpa\Pbll primal cauchyl

\Blem
\label{PreF}
Let $u\in\UR$.
With the above notations, the linearized system \Pbll\ has a unique solution $\soluzl\in\calY$.
Moreover, the estimate
\Beq
  \norma\soluzl_\calY \leq C \, \norma h_{\L2H}
  \label{preF}
\Eeq
holds true with a positive constant $C$ that depends only on the structure of the original system,
$\Omega$, $T$, the initial data, and \anold{$R$}.
\Elem

We do not give a \anold{detailed} proof of the lemma, \anold{restricting ourselves to 
sketch some formal estimates for brevity}.
We test equations \accorpa{primal}{terzal} by $\eta$, $\dt\psi$, and~$\zeta$, respectively, 
sum up, rearrange, and notice a cancellation.
Moreover, we add to both sides the same quantity 
$\frac12\frac d{dt}\iO|\psi|^2=\iO\psi\dt\psi$.
Then, the leading \lhs\ we obtain is given~by
\Beq
  \iO |\nabla\eta|^2
  + \tau \iO |\dt\psi|^2
  + \frac 12 \, \frac d{dt} \, \normaV\psi^2
  + \frac 12 \, \frac d{dt} \iO |\zeta|^2
  + \iO |\nabla\zeta|^2 \,.
  \non
\Eeq
In order to handle the products $\anold{\lambda_1}\psi\eta$ and $\anold{\lambda_2}\zeta\eta$ that enter the \rhs,
we argue as in Section~\ref{CONTDEP}.
Namely, we compute the mean value $\etaO$ of $\eta$ by taking $v=1$ in~\eqref{secondal} \anold{an arguing as done for \eqref{mean}. Namely, we square both sides, multiply the resulting equality by a suitable constant},
and repeatedly \juerg{invoke} the Young and Poincar\'e inequalities. 
Then, the \rhs\ can be easily dealt with, and estimate \eqref{preF} follows.
\anold{With this estimate at hand, the time derivative of $\zeta$ can be controlled in the dual space of $V$ by a comparison argument in \eqref{terzal}.}

\juerg{Next, we show that the linearized system analyzed above captures the Fr\'echet derivative of $\calS$.}

\Bthm
\label{Frechet}
Given any $u\in\UR$,
the \anold{solution operator} $\calS$ is \frechet\ differentiable at~$u$ 
\juerg{as a mapping from $\gianni{\UR\subset\calU}$  into~$\calY$}, 
and its \frechet\ derivative \juerg{$D\calS(u) \in \calL(\calU,\calY)$ associates}  
to every $h\in\calU$ the solution $\soluzl$ to the linearized system
corresponding to $u$ and to the variation~$h$.
\Ethm

\Bdim
The fact that the linear map described in the statement belongs to $\calL(\calU,\calY)$ 
is a consequence of Lemma~\ref{PreF}.
For the remainder of the proof, we assume without loss of generality that
$\norma h_\infty$ is small enough,
namely, such that the perturbed control $u+h$ also belongs to the open set~$\UR$,
so that the uniform bounds given by Theorem~\ref{Wellposedness}
also hold for the solution $\calS(u+\badtheta h)$ corresponding to $u+\badtheta h$
for every $\badtheta\in[0,1]$.
We set, for convenience,
\begin{align}
  & \anold{\soluzh := \calS(u+h), \quad
    \soluz := \calS(u)},
  \label{states}
  \\
  & \rho := \phih - \phi - \psi \,, \quad
  \theta := \muh - \mu - \eta,
  \aand
  \omega := \sigmah - \sigma - \zeta
  \label{rhothetaomega}
\end{align}
and observe that the triplet $(\rho,\theta,\omega)\in\calY$ solves the system
(all the equations holding \aet\ and for every $v\in V$)
\begin{align}
  & \iO \dt\rho \, v
  + \iO \nabla\theta \cdot \nabla v
  + m \iO \rho v
  = \iO \Lambda v\,,
  \label{primadiff}
  \\
  & \tau \iO \dt\rho \, v
  + \iO \nabla\rho \cdot \nabla v
  + \iO \Phi v
  - \iO \omega v
  = \iO \theta v\,,
  \label{secondadiff}
  \\
  &
	\anold{\<\dt\omega , v>_V}
  + \iO \nabla\omega \cdot \nabla v
  - \iO \Psi \cdot \nabla v
  = \iO \omega v
  - \iO \Xi \, v\,,
  \label{terzadiff}
  \\
  & \rho(0) = 0
  \aand
  \omega(0) = 0\,,
  \label{cauchydiff}
\end{align}
with the notation \eqref{deflambda} and \anold{with}
\begin{alignat}{2}
  & \Lambda := \gamma(\phih,\sigmah) - \gamma(\phi,\sigma) - \lamuno\psi - \lamdue\zeta \,, \quad
  && \Phi := F'(\phih)-F'(\phi)-\lambda\psi\,,
  \non
  \\
  & \Psi := \sigmah\nabla\phih - \sigma\nabla\phi - \zeta\nabla\phi - \sigma\nabla\psi \,, \quad
  && \Xi := \sigmah^2 - \sigma^2 - 2\sigma\zeta \,.
  \non
\end{alignat}
We \anold{then} test the equations \accorpa{primadiff}{terzadiff} by $\theta$, $\dt\rho$, and~$\omega$, respectively,
sum up, rearrange, and notice a cancellation.
Moreover, we add the same quantity $\frac12\frac d{dt}\iO|\rho|^2=\iO\rho\dt\rho$ to both sides.
We obtain, \aet,~that
\begin{align}
  & \iO |\nabla\theta|^2
  + \tau \iO |\dt\rho|^2
  + \frac 12 \, \frac d{dt} \, \normaV\rho^2
  + \frac 12 \, \frac d{dt} \iO |\omega|^2 
  + \iO |\nabla\omega|^2
  \non
  \\
  & = - m \iO \rho\theta
  + \iO \Lambda \, \theta
  - \iO \Phi \, \dt\rho
  + \iO \omega \dt\rho
  \non
  \\
  & \quad {}
  + \iO \Psi \cdot \nabla\omega
  + \iO |\omega|^2 
  - \iO \Xi \, \omega 
  + \iO \rho \dt\rho \,.
  \label{testeddiff}
\end{align}
At the same time, by arguing as in Section~\ref{CONTDEP},
we compute the mean value $\thetaO$ of $\theta$ \anold{from}~\eqref{secondadiff},
square the \juerg{resulting} equality, and divide both sides by a suitable constant.
Here, we have to manage the function $\Phi$ defined above, in addition,
since taking $v=1$ in~\eqref{secondadiff} yields
\Beq
  |\Omega| \, \thetaO 
  = \tau \iO \dt\rho
  + \iO \Phi
  - \iO \omega\,.
  \non
\Eeq 
By owing to the \anold{first-order} Taylor expansion of $F'$ with the \anold{second-order} remainder in integral form
(\anold{where} $F'''$ enters the remainder),
we have~that
\Beq
  \Phi = \lambda\rho + \Phi_0 
  \quad \hbox{with} \quad
  |\Phi_0| \leq c \, |\phih-\phi|^2 \,.
  \non
\Eeq
Therefore, proceeding as sketched above, we obtain~that
\Beq
  \frac {|\Omega|}{6\tau} \, |\thetaO|^2
  \leq \frac\tau 2 \iO |\dt\rho|^2
  + c \iO |\rho|^2
  + c \iO |\phih-\phi|^4
  + c \iO |\omega|^2\,, 
  \label{dasommare}
\Eeq
and we \anold{add} this inequality to \eqref{testeddiff}.
\juerg{We have the advantage} that we can use Poincar\'e's inequality
\juerg{to estimate} the terms on the \rhs\ of \eqref{testeddiff} involving~$\theta$.
Therefore, by repeatedly using Young's inequality, integrating with respect to time,
and applying the Gronwall lemma,
we can \anold{close the estimate}, provided we can suitably treat 
even the terms involving $\Lambda$, $\Phi$, $\Psi$\anold{,} and~$\Xi$.
To this end, we observe~that \anold{these can be rewritten as follows}:
\begin{align}
  & \Lambda = \lamuno\rho + \lamdue\omega + \Lambda_0
  \quad \hbox{with} \quad
  |\Lambda_0| \leq c \, \bigl( |\phih-\phi|^2 + |\sigmah-\sigma|^2 \bigr)\,,
  \non
  \\
  & \Psi = \omega \nabla\phi + (\sigmah-\sigma) \nabla(\phih-\phi) + \sigma\nabla\rho\,,
  \non
  \\
  & \Xi = |\sigmah-\sigma|^2 + 2\sigma\omega \,.
  \non
\end{align}
The expression \juerg{for} $\Lambda$ is just a \anold{first-order} Taylor expansion
(with the \anold{second-order} remainder in integral form, as above),
and the other\anold{s follow} \juerg{from} computing both sides
after eliminating $\phih$ and $\sigmah$ \juerg{using}~\eqref{rhothetaomega}.
Besides, we apply the second of the embeddings \anold{in}~\eqref{emb} and Theorem~\ref{Contdep} \anold{to infer that}
\begin{align}
  & \norma{\phih-\phi}_{\LQ4} 
  \leq c \, \norma{\phih-\phi}_{\L\infty H\cap\L2V}
  \leq c \, \norma h_{\L2H} \,,
  \non
  \\
  & \norma{\nabla(\phih-\phi)}_{\LQ4} 
  \leq c \, \norma{\phih-\phi}_{\L\infty V\cap\L2W}
  \leq c \, \norma h_{\L2H} \,,
  \non
  \\
  & \norma{\sigmah-\sigma}_{\LQ4} 
  \leq c \, \norma{\sigmah-\sigma}_{\L\infty H\cap\L2V}
  \leq c \, \norma h_{\L2H} \,.
  \non
\end{align}
Therefore, as it \anold{appears to be convenient to deal with the time integrated terms, we integrate the above estimate over time.}
By recalling that $\nabla\phi$ and $\sigma$ are bounded and that $\sigma$ is nonnegative \anold{as a consequence of Theorem~\ref{Wellposedness}}, 
we deduce that
\begin{align}
  & \intQt \Lambda \theta
  \leq \delta \, \iot |\thetaO(s)|^2 \, ds
  + \delta \intQt |\nabla\theta|^2
  \non
  \\
  & \quad {}
  + \cd \intQt \bigl( |\rho|^2+|\omega|^2 \bigr)
  + \cd \intQt \bigl( |\phih-\phi|^4 + |\sigmah-\sigma|^4 \bigr)
  \non
  \\
  \separa
  & \leq \delta \, \iot |\thetaO(s)|^2 \, ds
  + \delta \intQt |\nabla\theta|^2
  + \cd \intQt \bigl( |\rho|^2+|\omega|^2 \bigr)
  + \cd \, \norma h_{\L2H}^4\anold{\,,}
  \non
  \\
  & \quad 
  - \intQt \Phi \, \dt\rho
  \leq \delta \anold{\intQt}  |\dt\rho|^2
  + \cd \intQt |\rho|^2
  + \cd \intQt |\phih-\phi|^4
  \non
  \\
  \separa
  & \leq \delta \anold{\intQt} |\dt\rho|^2
  + \cd \intQt |\rho|^2
  + \cd \, \norma h_{\L2H}^4\anold{\,,}
  \non
  \\
  & \intQt \Psi \cdot \nabla\omega
  \leq \norma{\sigmah-\sigma}_{\LQ4} \, \norma{\nabla(\phih-\phi)}_{\LQ4} \, \norma{\nabla\omega}_{L^2(Q_t)}
  \leq \delta \intQt |\nabla\omega|^2
  + \cd \, \norma h_{\L2H}^4\anold{\,,}
  \non
  \\
  & \quad 
  - \intQt \Xi \, \omega
  \leq - \intQt |\sigmah-\sigma|^2 \omega
  \leq \intQt |\omega|^2 
  + \anold{c}\intQt |\sigmah-\sigma|^4
  \leq \intQt |\omega|^2 
  + c \, \norma h_{\L2H}^4 \,.
  \non
\end{align}
\anold{Upon collecting the above estimates,}  choosing $\delta$ small enough\anold{,} and applying the Gronwall lemma\anold{,} we conclude~that
\begin{align}
  & \norma\rho_{\H1H\cap\L\infty V}^2
  + \norma\theta_{\L2V}^2
  + \norma\omega_{\L\infty H\cap\L2V}^2
  \non
  \\
  & \leq c \, \norma h_{\L2H}^4 
  \leq c \, \norma h_\calU^4 \,.
  \non
\end{align}
\anold{Due to the above estimate, it is then a standard matter to infer from a comparison argument in \eqref{terzadiff} that 
\begin{align*}
	\norma{\dt \omega}_{\L2 \Vp}^2
 \leq  c \, \norma h_{\L2H}^4 
  \leq c \, \norma h_\calU^4 \,.
\end{align*}}
\anold{Thus, we conclude} that
\Beq
  \lim_{\norma h_\calU\to0} \frac{\norma{\calS(u+h)-\calS(u)-\soluzl}_\calY}{\norma h_\calU}
  = 0\,,
  \non
\Eeq
and this is the assertion of the statement.
\Edim


\subsection{First-order optimality conditions}
\label{OPT}

The above \frechet\ differentiability result
permits us to apply the chain rule to the composite map
\Beq
  \UR \ni u \mapsto \anold{(u,\phi,\sigma)} \mapsto \calJ(u\anold{;}\phi,\sigma) \anold{\,,}
  \non
\Eeq
\anold{where $\phi$ and $\sigma$ are components of the solution to the state system \Ipbl\ corresponding to the control variable~$u$.}
Since $\Uad$ is convex,
one immediately sees that a necessary condition for 
$\ustar$ to be an optimal control is given by the variational inequality
\begin{align}
  & \alpha_1 \intQ (\phistar-\tarQ) \psi
  + \alpha_2 \iO (\phistar(T)-\tarO) \psi(T)
    + \alpha_3 \intQ (\sigmastar-\tarsQ) \zeta
	\non\\	
  &\quad+ \alpha_4 \iO (\sigmastar(T)-\tarsO) \zeta(T)
    + \alpha_5 \intQ \ustar (u-\ustar)\,\ge\,0 \,\quad \hbox{for every $u\in\Uad$}\,,
  \label{preopt}
\end{align}
where $\soluzl$ is the solution to the linearized system \Pbll\ associated 
with the \anold{state} $\soluzstar:=\calS(\ustar)$ and the \anold{increment} $h=u-\ustar$.

\anold{Unfortunately, the above variational inequality is not} \juerg{helpful, since it requires to solve
the linearized problem infinitely many times. As usual, this difficulty is bypassed by introducing a 
proper adjoint problem}. To this end,
\anold{let us} fix an optimal control $\ustar\in\Uad$ and its corresponding state $\soluzstar=\calS(\ustar)$
and upgrade the abbreviations \eqref{deflambda} by setting
\Beq
  \lamuno := \gamma_\phi(\phistar,\sigmastar) , \quad
  \lamdue := \gamma_\sigma(\phistar,\sigmastar),
  \aand
  \lam := F''(\phistar) \,.
  \label{deflambdastar}
\Eeq 
Furthermore, we define
\begin{align}
  & g_1 := \alpha_1(\phistar-\tarQ) , \quad
  g_2 := \alpha_2(\sigmastar-\tarsQ) , \quad
  g_3 := \alpha_3 \bigl( \phistar(T)-\tarO \bigr),
  \non
  \\
  & \aand
  g_4 := \alpha_4 \bigl( \sigmastar(T)-\tarsO \bigr) \,.
  \label{defg}
\end{align}
Then, the associated adjoint problem consists in finding a triplet $\soluza$
satisfying
\begin{align}
  & p \in \L2V , \quad 
  q \in \L2V , \quad
  p+\tau q\in\H1\Vp,
  \non
  \\
  & \aand
  r \in \H1\Vp\anold{\cap \L\infty H}\cap\L2V,
  \label{regsoluza}
\end{align}
\juerg{which solves the} backward-in-time system
\begin{align}
  & - \< \dt(p+\tau q) , v >_V
  + \iO \nabla q \cdot \nabla v
  + \iO \bigl( (m-\lamuno)p + \lam q \bigr) v
  - \iO \sigmastar \nabla r \cdot \nabla v
  \non
  \\
  & = \iO g_1 v
  \quad \hbox{\aet, for every $v\in V$}\,,
  \label{primaa}
  \\
  & \iO \nabla p \cdot \nabla v
  = \iO q v
  \quad \hbox{\aet, for every $v\in V$}\,,
  \label{secondaa}
  \\
  & - \< \dt r , v >_V
  + \iO \nabla r \cdot \nabla v
  + \iO \bigl( -\lamdue p - q + (2\sigmastar-1) r \bigr) v
  - \iO \nabla\phistar \cdot \nabla r \, v
  \non
  \\
  & = \iO g_2 v
  \quad \hbox{\aet, for every $v\in V$}\,,
  \label{terzaa}
  \\
  & (p+\tau q)(T) = g_3
  \aand
  r(T) = g_4 \,.
  \label{cauchya}
\end{align}
\Accorpa\Pbla primaa cauchya

\Bthm
\label{Wellposednessa}
The adjoint problem \Pbla\ has a unique solution $\soluza$ satisfying 
the regularity requirements \anold{in} \eqref{regsoluza}.
\Ethm

\Bdim
To \anold{show} well-posedness, it is convenient to introduce the auxiliary unknown 
\Beq
  z := p+\tau q
  \label{defz}
\Eeq
and \juerg{to} eliminate $q$ from the system. \juerg{We are going to prove the well-posedness of \Pbla\ written in 
terms of the new variables, which is equivalent.}
With the above transformation, the new \juerg{sought} triplet is $(z,p,r)$, which is required to satisfy
\begin{align}
  & z \in \H1\Vp \anold{\cap \L\infty H}\cap \L2V , \quad
  p \in \L2V,
  \non
  \\
  & \aand
  r \in \H1\Vp \anold{\cap \L\infty H}\cap \L2V.
  \label{regsoluzz}
\end{align}
The adjoint problem \Pbla\ \juerg{then takes} the form
\begin{align}
  & - \< \dt z , v >_V
  + \frac 1\tau \iO \nabla z \cdot \nabla v
  - \frac 1\tau \iO \nabla p \cdot \nabla v
  \non
  \\
  & \quad {}
  + \iO \Bigl( (m-\lamuno)p + \lam \, \frac{z-p}\tau \Bigr) v
  - \iO \sigmastar \nabla r \cdot \nabla v
  \non
  \\
  & = \iO g_1 v
  \quad \hbox{\aet, for every $v\in V$}\,,
  \label{primaz}
  \\
  \separa
  & \iO \nabla p \cdot \nabla v
  + \frac 1\tau \iO p v
  = \frac 1\tau \iO z v
  \quad \hbox{\aet, for every $v\in V$}\,,
  \label{secondaz}
  \\
  & - \< \dt r , v >_V
  + \iO \nabla r \cdot \nabla v
  + \iO \Bigl( -\lamdue p - \frac{z-p}\tau + (2\sigmastar-1) r \Bigr) v
  - \iO \nabla\phistar \cdot \nabla r \, v
  \non
  \\
  & = \iO g_2 v
  \quad \hbox{\aet, for every $v\in V$}\,,
  \label{terzaz}
  \\
  \separa
  & z(T) = g_3
  \aand
  r(T) = g_4 \,.
  \label{cauchyz}
\end{align}
\Accorpa\Pblz primaz cauchyz

\juerg{We now prove the well-posedness for the new problem}.
\juerg{To show existence}, we start from a regularized version obtained as follows\anold{: we} 
divide \eqref{secondaz} by $\tau$ and add a viscosity term \anold{$-\eps \dt p^\eps$}
depending on the parameter $\eps>0$ on the \lhs.
Moreover, we multiply \eqref{terzaz} by a constant $N$ whose value will be chosen later on.
Finally, we rewrite the three equations with independent test functions $v_1$, $v_2$, and $v_3$,
and add them to each other. 
\juerg{This leads us to the following approximating problem: find} a triplet $\soluzeps$,
\juerg{which satisfies}
\Beq
  \zeps,\, \peps,\, \reps \in \H1\Vp \anold{\cap \L\infty H}\cap \L2V
  \label{regeps}
\Eeq
\juerg{and solves the variational equation}
\begin{align}
  & - \< \dt\zeps , v_1 >_V
  + \frac 1\tau \iO \nabla\zeps \cdot \nabla v_1
  - \frac 1\tau \iO \nabla\peps \cdot \nabla v_1
  \non
  \\
  & \quad {}
  + \iO \Bigl( (m-\lamuno)\peps + \lam \, \frac{\zeps-\peps}\tau \Bigr) v_1
  - \iO \sigmastar \nabla\reps \cdot \nabla v_1
  \non
  \\
  & \quad {}
  - \eps \< \dt\peps , v_2 >_V
  + \frac 1\tau\iO \nabla\peps \cdot \nabla v_2
  + \frac 1{\tau^2} \iO \peps v_2
  - \frac 1{\tau^2} \iO \zeps v_2
  \non
  \\
  & \quad {}
  - N \< \dt\reps , v_3 >_V
  + N \iO \nabla\reps \cdot \nabla v_3
  \non
  \\
  & \quad {}
  + N \iO \Bigl( -\lamdue\peps - \frac{\zeps-\peps}\tau + (2\sigmastar-1) \reps \Bigr) v_3
  - N \iO \nabla\phistar \cdot \nabla\reps \, v_3
  \non
  \\  
  & = 
  \iO g_1 v_1
  + N \iO g_2 v_3
  \qquad \hbox{\aet, for every $v_1,\,v_2,\,v_3\in V$}\,,
  \label{eqeps}
  \\[2mm]
  & (\zeps,\peps,\reps)(T)
  = (g_3,0,g_4) \,.
  \label{cauchyeps}
\end{align}
\Accorpa\Pbleps eqeps cauchyeps
In order to see that this problem is well\anold{-}posed,
we introduce the Hilbert triplet  
\Beq
  (\calV,\calH,\calVp),
  \quad \hbox{where} \quad
  \calV := V\juerg{\times V\times V}
  \aand
  \calH := H\juerg{\times H\times H},
  \label{calVH}
\Eeq
with the embedding $\calH\emb\calVp$ that is associated 
with the \juerg{following inner product in~$\calH$ (which is equivalent to the standard one):}
\begin{align}
  & (w,v)_\calH
  := \iO w_1 v_1
  + \eps \iO w_2 v_2
  + N \iO w_3 v_3
  \non
  \\
  & \quad \hbox{for $w=(w_1,w_2,w_3),\ v=(v_1,v_2,v_3)\in\calH$}.
  \label{prodH}
\end{align}
Notice that, for $w=(w_1,w_2,w_3)\in\calH$, we have~that
\Beq
  \< w,v >_\calV
  = \< w_1,v_1>_V + \eps \< w_2,v_2 >_V + N \< w_3,v_3 >_V
  \quad \hbox{for every $v=(v_1,v_2,v_3)\in\calV$},
  \non
\Eeq
since $\<w,v>_\calV=(w,v)_\calH$ and $\<w_i,v_i>_V=\iO w_iv_i$\anold{, $i=1,2,3$}.
Therefore, the same relation holds true for every $w=(w_1,w_2,w_3)\in\calVp=\juerg{\Vp\times\Vp\times\Vp}$,
so that \eqref{eqeps} takes the~form
\begin{align}
  & - \< \dt\soluzeps , v >_\calV
  + a \bigl( \cpto\anold{;}  \soluzeps , v \bigr)
  + b \bigl( \cpto\anold{;} \soluzeps , v \bigr)
  \non
  \\
  & = \bigl( (\juerg{g_1},0,\juerg{g_2}) , v \bigr)_\calH
  \quad \hbox{\aet, for every $v\in\calV$}\,,
  \label{astratta}
\end{align}
where $a$ and $b$ are the time-dependent continuous bilinear forms on $\calV\times\calV$ defined~by
\begin{align}
  & a(t;w,v)
  := \frac 1\tau \iO \nabla w_1 \cdot \nabla v_1
  - \frac 1\tau \iO \nabla w_2 \cdot \nabla v_1
  - \iO \sigmastar(t) \nabla w_3 \cdot \nabla v_1
  \non
  \\
  & \quad{}
  + \frac 1\tau \iO \nabla w_2 \cdot \nabla v_2
  + N \iO \nabla w_3 \cdot \nabla v_3
  - N \iO \nabla\phistar(t) \cdot \nabla w_3 \, v_3
  \label{defa}
  \\
  & b(t;w,v)
  := \iO \Bigl( (m-\lamuno(t)) w_2 + \lam(t) \, \frac{w_1-w_2}\tau \Bigr) v_1
  + \frac 1{\tau^2} \iO w_2 v_2
  - \frac 1{\tau^2} \iO w_3 v_2
  \non
  \\
  & \quad
  + N \iO \Bigl( -\lamdue(t) w_2 - \frac{w_1-w_2}\tau + (2\sigmastar(t)-1) w_3 \Bigr) v_3
  \label{defb}
\end{align}
both \aat\ and $w=(w_1,w_2,w_3),\ v=(v_1,v_2,v_3)\in\calV$.
We notice at once that
\begin{align}
  & |a(t;w,v)|
  \leq c \, \norma w_\calV \, \norma v_\calV
  \aand
  |b(t;w,v)|
  \leq c \, \norma w_\calH \, \norma v_\calH
  \non
  \\
  & \quad \hbox{\aat\ and every $w,v\in\calV$}\,,
  \label{bbnd}
\end{align}
since $\lamuno$, $\lamdue$, $\lam$ and $\sigmastar$ are bounded.
Moreover, $(g_1,0,g_2)\in\L2\calH$ and $(g_3,0,g_4)\in\calH$.
Therefore, the existence of a unique solution 
\Beq
  \soluzeps \in \H1\calVp \anold{\cap \L\infty \calH}\cap \L2\calV
  \non
\Eeq
is ensured if we can find positive constants $\alpha$ and $\kappa$ such that
\Beq
  a(t;w,w) + \kappa \, \norma w_\calH^2
  \geq \alpha \, \norma w_\calV^2 
  \quad \hbox{\aat\ and every $\anold{w}\in\calV$} \,.
  \label{coercive}
\Eeq
To prove this, we account for the Young inequality,
(which implies, in particular, that $|x|^2-x\cdot y+|y|^2\geq\frac12(|x|^2+|y|^2)$ for every $x,y\in\erre^2$) 
and suitably choose the value of~$N$.
We have that
\begin{align}
  & a(t;w,w)
  = \frac 1\tau \iO |\nabla w_1|^2
  - \frac 1\tau \iO \nabla w_2 \cdot \nabla w_1
  - \iO \sigmastar(t) \nabla w_3 \cdot \nabla w_1
  \non
  \\
  & \quad{}
  + \frac 1\tau \iO |\nabla w_2|^2 
  + N \iO |\nabla w_3|^2
  - N \iO \nabla\phistar(t) \cdot \nabla w_3 \, w_3
  \non
  \\
  & \geq \frac 1{2\tau} \iO \bigl( |\nabla w_1|^2 + |\nabla w_2|^2 \bigr)
  + N \iO |\nabla w_3|^2
  \non
  \\
  & \quad {}
  - \frac 1{4\tau} \iO |\nabla w_1|^2
  - \tau \, \norma\sigmastar_\infty^2 \, \iO |\nabla w_3|^2
  - \frac N2 \iO |\nabla w_3|^2 
  - \frac N2 \, \norma{\nabla\phistar}_\infty^2 \, \iO |w_3|^2 \,.
  \label{percoerc}
\end{align}
Therefore, \eqref{coercive} is satisfied if we choose
\Beq
  N = 2\tau \, \norma\sigmastar_\infty^2 + 1 , \quad
  \alpha = \min\{ 1/(4\tau),1/2 \}
  \aand
  \kappa = \frac N2 \, \norma{\nabla\phistar}_\infty^2\anold{.}
  \label{defNak}
\Eeq

\def\intQt{\int_{Q^t}}%

\anold{The next step consists of showing some uniform (with respect to~$\eps$) 
estimates that allow us to let the approximating parameter tend to zero.}
We test \eqref{astratta} by $\soluzeps$ and integrate over~$(t,T)$\anold{, for an arbitrary $t \in [0,T)$}.
We obtain an equality and notice that the unique term that contains~$\eps$ is nonnegative.
Therefore, by ignoring it, we deduce that
(recall \eqref{defQt} for the definition of~$Q^t$)
\begin{align}
  & \frac 12 \iO |\zeps(t)\anold{|}^2
  + \frac N2 \iO |\reps(t)|^2
  \non
  \\
  & \quad {}
  + \int_t^T a(s;\soluzeps(s),\soluzeps(s)) \, ds
  + \int_t^T b(s;\soluzeps(s),\soluzeps(s)) \, ds
  \non
  \\
  & \leq \frac 12 \iO |g_3|^2
  + \frac N2 \iO |g_4|^2
  + \intQt g_1 \, \zeps
  + \intQt g_2 \, \reps \,,
  \non
\end{align}
\juerg{where $\eps$ no longer appears}.
However, since a leading term is now missing, we cannot simply apply \eqref{bbnd} and \eqref{coercive}.
Nevertheless, by coming back to the inequality \anold{\eqref{percoerc}}, 
we see~that
\Beq
  \int_t^T a(s;\soluzeps(s),\soluzeps(s)) \, ds
  \geq \alpha \intQt \bigl( |\nabla\zeps|^2 + |\nabla\peps|^2 + |\nabla\reps|^2 \bigr)
  - \kappa \intQt |\reps|^2
  \non
\Eeq
with $\alpha$ and~$\kappa$ given \anold{as in} \eqref{defNak}.
Moreover, we observe that the integral involving $b$ contains the quantity
$(1/\tau^2)|\peps|^2$, and this can replace the missing leading term, 
since it yields a positive contribution to the \lhs.
It follows that the remaining terms involving $\peps$ can be dealt with 
using Young's inequality.
Therefore, applying the (backward) Gronwall lemma, 
we conclude~that
\Beq
  \norma\zeps_{\L\infty H\cap\L2V}
  + \norma\peps_{\L2V}
  + \norma\reps_{\L\infty H\cap\L2V}
  \leq c \,.
  \label{stimaadj}
\Eeq
At this point, we come back to \eqref{eqeps},
test it by $(v_1,0,v_3)$ with arbitrary time-dependent test functions $v_1,v_3\in\L2V$
and easily deduce from \eqref{stimaadj}~that
\Beq
  \norma{\dt\zeps}_{\L2\Vp}
  + \norma{\dt\reps}_{\L2\Vp}
  \leq c \,.
  \non
\Eeq
By this and \eqref{stimaadj}, and owing to well-known compactness results,
we see that, at least for a \anold{sub}sequence $\{\eps_n\}\seto0$
(however, we still write $\eps$ for simplicity),
we have~that
\begin{align}
  & \zeps \to z
  \quad \hbox{weakly star in $\H1\Vp\anold{\cap \L\infty H}\cap\L2V$}\,,
  \label{convz}
  \\
  & \peps \to p
  \quad \hbox{weakly in $\L2V$}\,,
  \label{convp}
  \\
  & \reps \to r
  \quad \hbox{weakly star in $\H1\Vp\anold{\cap \L\infty H}\cap\L2V$}\,,
  \label{convr}
\end{align}
for some triplet $\soluzaz$ satisfying \eqref{regsoluzz}.
Then, it is clear that this triplet also satisfies \eqref{primaz} and~\eqref{terzaz}
(since all of the coefficients are bounded), as well as \eqref{cauchyz}.
To prove that \eqref{secondaz} is fulfilled, 
we take any $v\in H^1_0(0,T;V)$, test \eqref{astratta} by $(0,\tau v,0)$,
and integrate over~$(t,T)$.
We obtain~that
\Beq
  - \tau\eps \ioT \< \dt\peps(t) , v(t) >_V \, dt
  + \intQ \nabla\peps \cdot \nabla v
  + \frac 1\tau \intQ \peps v
  - \frac 1\tau \intQ \anold{\zeps} v
  = 0 \,.
  \non
\Eeq
\anold{Moreover, s}ince
\Beq
  - \tau\eps \ioT \< \dt\peps(t) , v(t) >_V \, dt
  = \tau\eps \ioT \< \peps(t) , \dt v(t) >_V \, dt
  = \tau\eps \intQ \peps \dt v\,,
  \non
\Eeq
we can let $\eps$ tend to zero on account of \eqref{convp} and~\eqref{convr}.
\anold{Finally, we deduce}~that
\Beq
  \intQ \nabla p \cdot \nabla v
  + \frac 1\tau \intQ p v
  - \frac 1\tau \intQ \anold{z} v
  = 0 \,\anold{,}
  \non
\Eeq
\anold{which, s}ince $v\in H^1_0(0,T;V)$ is arbitrary, this is equivalent to~\eqref{secondaz}.

\juerg{Finally}, we prove uniqueness for \Pblz.
%
Since the problem is linear, \an{it suffices to prove uniqueness for the homogeneous case.
Hence,} we replace the known terms $g_i$ by zero\an{, $i=1,...,4$.}
Then, we test the equations by $z$, $p$ and~$Nr$, respectively,
where $N$ is given by~\eqref{defNak}.
After some rearragement\anold{s}, we obtain~that
\begin{align}
  & - \frac 12 \, \frac d{dt} \iO |z|^2
  - \frac N2 \, \frac d{dt} \iO |r|^2
  \non
  \\
  & \quad {}
  + \frac 1\tau \iO |\nabla z|^2
  - \frac 1\tau \iO \nabla p \cdot \nabla z
  - \iO \sigmastar \nabla r \cdot \nabla z
  + \iO |\nabla p|^2
  + N \iO |\nabla r|^2
  \non
  \\
  & \quad {}
  + \frac 1\tau \iO |p|^2
  \non
  \\
  & = \frac 1\tau \iO z p
  - \iO \Bigl( (m-\lamuno)p + \lam \, \frac{z-p}\tau \Bigr) z
  \non
  \\
  & \quad {}
  - N \iO \Bigl( -\lamdue p - \frac{z-p}\tau + (2\sigmastar-1) r \Bigr) r
  + N \iO \nabla\phistar \cdot \nabla r \, r \,.
  \non
\end{align}
Now, we recall \eqref{percoerc} and apply it to $\soluzaz$ 
to estimate the second line of the above equality from below.
At the same time, we repeatedly account for the Young inequality
to estimate the \rhs.
Then, we integrate over $(t,T)$ and apply the (backward) Gronwall lemma
\anold{which} yields $\soluzaz=(0,0,0)$,
and the proof is complete.
\Edim

We are now ready to prove a satisfactory version of the \anold{first-order} necessary condition for optimality.

\Bthm
\label{Opt}
Let $\ustar$ and $\soluzstar$ be an optimal control and the corresponding state, respectively.
Moreover, with the notations \accorpa{deflambdastar}{defg},
let $\soluza$ be the unique solution to the corresponding adjoint problem \Pbla.
Then, there holds the variational inequality
\Beq
  \intQ (r+\alpha_5\ustar) (u-\ustar) \geq 0
  \quad \hbox{for every $u\in\Uad$}.
  \label{opt}
\Eeq
In particular, if $\alpha_5>0$, then $\ustar$ is the $\LQ2$-projection of $-r/\alpha_5$ on~$\Uad$.
\Ethm

\Bdim
To prove \eqref{opt}, we fix $u\in\Uad$ 
and consider the linearized system \Pbll\ associated with $\soluz=\soluzstar$
(thus with $\lambda_1$, $\lambda_2$ and $\lambda$ given by~\eqref{deflambdastar})
\juerg{and} $h=u-\ustar$.
By testing the equations by $p$, $q$ and~$r$, respectively, and summing up,
we obtain~that
\begin{align}
  & \iO \dt\psi \, p
  + \iO \nabla\eta \cdot \nabla p
  + m \iO \psi p
  \non
  \\
  & \quad {}
  + \tau \iO \dt\psi \, q
  + \iO \nabla\psi \cdot \nabla q
  + \iO \bigl( \lam \psi - \zeta \bigr) q
  \non
  \\
  & \quad {}
  \anold{+ \<\dt\zeta , r>_V}
  + \iO \nabla\zeta \cdot \nabla r
  - \iO \bigl( \zeta\nabla\phistar + \sigmastar\nabla\psi \bigr) \cdot \nabla r
  \non
  \\
  \separa
  & = \iO \bigl( \lamuno \psi + \lamdue \zeta \bigr) p
  + \iO \eta q
  + \iO \bigl( \zeta - 2\sigmastar\zeta + (u-\ustar) \bigr) r\anold{.}
  \label{testlin}
\end{align}
At the same time, we test the equations \accorpa{primaa}{terzaa} of the adjoint system
by $-\psi$, $-\eta$, and~$-\zeta$, respectively,
and add them to each other.
We \anold{find}~that
\begin{align}
  & \< \dt(p+\tau q) , \psi >_V
  - \iO \nabla q \cdot \nabla\psi
  - \iO \bigl( (m-\lamuno)p + \lam q \bigr) \psi
  + \iO \sigmastar \nabla r \cdot \nabla\psi
  \non
  \\
  & \quad {}
  - \iO \nabla p \cdot \nabla \eta
  + \< \dt r , \zeta >_V
  - \iO \nabla r \cdot \nabla\zeta
  \non
  \\
  & \quad {}
  - \iO \bigl( -\lamdue p - q + (2\sigmastar-1) r \bigr) \zeta
  + \iO \nabla\phistar \cdot \nabla r \, \zeta
  \non
  \\
  & = - \iO q \eta
  - \iO g_1 \psi
  - \iO g_2 \zeta\anold{.}
  \label{testadj}
\end{align}
Now, we take the sum of \eqref{testlin} and \eqref{testadj}
\anold{and notice that some cancellations} occur\anold{. It then} remains~that
\begin{align}
  & \iO \dt\psi (p+\tau q)
  \anold{+ \<\dt\zeta , r>_V}
  + \< \dt(p+\tau q) , r >_V
  + \< \dt r , \zeta >_V
  \non
  \\
  & = \iO (u-\ustar) \, r
  - \iO g_1 \psi
  - \iO g_2 \zeta \,.
  \non
\end{align}
\anold{Thus}, we integrate over $(0,T)$ and \juerg{employ} the well-known integration-by-parts formula
for functions belonging to $\H1\Vp\cap\L2V$.
Recalling the Cauchy conditions \eqref{cauchyl} and~\eqref{cauchya},
we deduce~that
\Beq
  \iO \psi(T) \, g_3
  + \iO \zeta(t) \, g_4
  = \intQ (u-\ustar) \, r
  - \intQ g_1 \, \psi
  - \intQ g_2 \, \zeta \,.
  \non
\Eeq
By recalling the definition \eqref{defg} of the functions $g_i$
and \juerg{inserting the last formula in}~\eqref{preopt},
we obtain~\eqref{opt}\juerg{, which concludes the proof}.
\Edim


\section*{Acknowledgments}
\anold{
AS gratefully \juerg{acknowledges} partial support 
from the MIUR-PRIN Grant 2020F3NCPX ``Mathematics for industry 4.0 (Math4I4)''
and his affiliation to the GNAMPA (Gruppo Nazionale per l'Analisi Matematica, 
la Probabilit\`a e le loro Applicazioni) of INdAM (Isti\-tuto 
Nazionale di Alta Matematica).
}


\vspace{3truemm}


\vspace{3truemm}

\Begin{thebibliography}{10} \footnotesize

\bibitem{AC}
A. R. A. Anderson and M. A. J. Chaplain. 
Continuous and discrete mathematical models of tumor-induced angiogenesis. 
{\it Bull. Math. Biol.}, {\bf 60 (5)} (1998), 857–899.

\bibitem{ALL}
A. Agosti, A. G. Lucifero and S. Luzzi.
An image-informed Cahn--Hilliard Keller--Segel multiphase field model for tumor growth with angiogenesis. {\it Appl. Math. Comput.}, {\bf 445C} (2023), 127834.

\bibitem{AS}
A. Agosti and A.~Signori.
Multi-species Cahn–Hilliard–Keller–Segel tumor growth model with chemotaxis and angiogenesis.
{\it In preparation}.

\bibitem{CSS}
P. Colli, A. Signori and J. Sprekels.
Optimal control of a phase field system modelling tumor growth with chemotaxis and singular potentials.
{\it Appl. Math. Optim.} {\bf 83} (2021), 2017-2049.
%
\bibitem{CSS2}
P. Colli, A. Signori and J. Sprekels.
Second-order analysis of an optimal control problem in a phase field tumor growth model with singular potentials and chemotaxis.
{\it ESAIM Control Optim. Calc. Var.}, {\bf 27} (2021). doi.org/10.1051/cocv/2021072.	

\bibitem{CSS4}
P. Colli, A. Signori and J. Sprekels.
Optimal control problems with sparsity for phase field tumor growth models involving variational inequalities.
{\it J. Optim. Theory Appl.}, (2022). doi.org/10.1007/s10957-022-02000-7.

\bibitem{EG}
M.~Ebenbeck and H.~Garcke.
Analysis of a Cahn--Hilliard--Brinkman model for tumour growth with chemotaxis.
{\it J. Differ. Equ.}, {\bf 266 (9)} (2019),  5998-6036.

\bibitem{FK}
J. Folkman. 
Tumor angiogenesis. 
{\it Adv. Cancer Res.}, {\bf 43} (1985), 175–203.

\bibitem{GARL_2}
H. Garcke and K. F. Lam.
Analysis of a Cahn--Hilliard system with non--zero Dirichlet 
conditions modeling tumor growth with chemotaxis.
{\it Discrete Contin. Dyn. Syst.} {\bf 37 (8)} (2017), 4277-4308.

\bibitem{GLSS}
H.~Garcke, K. F.~Lam, E.~Sitka and V.~Styles.
A Cahn--Hilliard--Darcy model for tumour growth with chemotaxis and active transport.
{\it Math. Models Methods Appl. Sci.}, {\bf 26 (6)} (2016),  1095-1148.

\bibitem{H}
D. Horstmann.
From 1970 until now: the Keller--Segel model in chemotaxis and its consequences.
{\it Jahresber. Deutsch. Math.-Verei.}, {\bf 106} (2004), 51-69.

\bibitem{KS2}
P. Knopf and A. Signori.
Existence of weak solutions to multiphase Cahn--Hilliard--Darcy and Cahn--Hilliard--Brinkman models for stratified tumor growth with chemotaxis and general source terms. 
{\it Comm. Partial Differential Equations}, {\bf 47 (2)} (2022), 233-278.

\bibitem{RSS1}
E.~Rocca, L.~Scarpa and A.~Signori.
Parameter identification for nonlocal phase field models for tumor growth via optimal control and asymptotic analysis. 
{\it Math. Models Methods Appl. Sci.}, {\bf 31 (13)} (2021), 2643-2694. 

\bibitem{RSS}
E. Rocca, G. Schimperna and A. Signori.
On a Cahn--Hilliard--Keller--Segel model with generalized logistic source 
describing tumor growth.
{\it J.~Differ. Equ.}, {\bf 343} (2023)\an{,} 530-578.

\bibitem{SS}
L.~Scarpa and A.~Signori.
On a class of non-local phase-field models for tumor growth with possibly singular potentials, chemotaxis, and active transport.
{\it Nonlinearity},  {\bf 34} (2021), 3199-3250.

\bibitem{S}
A. Signori.
Optimal distributed control of an extended model of tumor
growth with logarithmic potential.
{\it Appl. Math. Optim.}, {\bf 82} (2020), 517-549.

\bibitem{Winkler1}
M. Winkler.
Boundedness in the higher-dimensional parabolic-parabolic chemotaxis system with logistic source.
{\it Comm. Partial Differential Equations}, {\bf 35} (2010), 1516–1537.

\bibitem{Winkler2}
M. Winkler.
Finite-time blow-up in the higher-dimensional parabolic–parabolic Keller--Segel system.
{\it J. Math. Pures Appl.}, {\bf 100} (2013), 748–767.

\bibitem{Winkler3}
M. Winkler.
Emergence of large population densities despite logistic growth restrictions in fully parabolic chemotaxis systems.
{\it Discrete Contin. Dyn. Syst. Ser. B}, {\bf 22} (2017), 2777–2793.

\End{thebibliography}

\End{document}


Moreover, $\soluza$ also enjoys the further regularity properties
\Beq
  p \in \L2{W\cap\Hx3}
  \aand
  r \in \L2W \,.
  \label{morereg}
\Eeq